\documentclass[english]{IEEEtran}
\usepackage[T1]{fontenc}
\usepackage[latin9]{inputenc}
\usepackage{color}
\usepackage{float}
\usepackage{amsmath}
\usepackage{amsthm}
\usepackage{amssymb}
\usepackage{graphicx}

\makeatletter

\providecommand{\tabularnewline}{\\}
\floatstyle{ruled}
\newfloat{algorithm}{tbp}{loa}
\providecommand{\algorithmname}{Algorithm}
\floatname{algorithm}{\protect\algorithmname}

\theoremstyle{definition}
\newtheorem{example}{\protect\examplename}
\theoremstyle{definition}
\newtheorem{defn}{\protect\definitionname}
\theoremstyle{plain}
\newtheorem{thm}{\protect\theoremname}
\theoremstyle{plain}
\newtheorem{lem}{\protect\lemmaname}

\usepackage{cite}

\newtheorem{assumption}{Assumption}

\usepackage{subfig}

\makeatother

\usepackage{babel}
\providecommand{\definitionname}{Definition}
\providecommand{\examplename}{Example}
\providecommand{\lemmaname}{Lemma}
\providecommand{\theoremname}{Theorem}

\begin{document}
\title{Two-Stage Stochastic Optimization via Primal-Dual Decomposition and
Deep Unrolling }
\author{{\normalsize{}An Liu, }\textit{\normalsize{}Senior Member, IEEE}{\normalsize{},
Rui Yang, Tony Q. S. Quek,}\textit{\normalsize{} Fellow, IEEE}{\normalsize{}
and Min-Jian Zhao, }\textit{\normalsize{}Member, IEEE}{\normalsize{}}\thanks{This work was supported in part by the National Science Foundation
of China under Grant 62071416, in part by the National Research Foundation,
Singapore and Infocomm Media Development Authority under its Future
Communications Research \& Development Programme, in part by the SUTD
Growth Plan Grant for AI, and in part by the SUTD-ZJU Seed Grant SUTD-ZJU
(SD) 201909. Any opinions, findings and conclusions or recommendations
expressed in this material are those of the author(s) and do not reflect
the views of National Research Foundation, Singapore and Infocomm
Media Development Authority. (Corresponding authors: An Liu; Min-Jian
Zhao.)

An Liu, Rui Yang and Min-Jian Zhao are with the College of Information
Science and Electronic Engineering, Zhejiang University, Hangzhou
310027, China (email: anliu@zju.edu.cn).

Tony Q. S. Quek is with the Information Systems Technology and Design
Pillar, Singapore University of Technology and Design (email: tonyquek@sutd.edu.sg).}}
\maketitle
\begin{abstract}
We consider a two-stage stochastic optimization problem, in which
a long-term optimization variable is coupled with a set of short-term
optimization variables in both objective and constraint functions.
Despite that two-stage stochastic optimization plays a critical role
in various engineering and scientific applications, there still lack
efficient algorithms, especially when the long-term and short-term
variables are coupled in the constraints. To overcome the challenge
caused by tightly coupled stochastic constraints, we first establish
a two-stage primal-dual decomposition (PDD) method to decompose the
two-stage problem into a long-term problem and a family of short-term
subproblems. Then we propose a PDD-based stochastic successive convex
approximation (PDD-SSCA) algorithmic framework to find KKT solutions
for two-stage stochastic optimization problems. At each iteration,
PDD-SSCA first runs a short-term sub-algorithm to find stationary
points of the short-term subproblems associated with a mini-batch
of the state samples. Then it constructs a convex surrogate for the
long-term problem based on the deep unrolling of the short-term sub-algorithm
and the back propagation method. Finally, the optimal solution of
the convex surrogate problem is solved to generate the next iterate.
We establish the almost sure convergence of PDD-SSCA and customize
the algorithmic framework to solve two important application problems.
Simulations show that PDD-SSCA can achieve superior performance over
existing solutions. 
\end{abstract}


\section{Introduction}

In this paper, we consider the following two-stage stochastic optimization
problem:
\begin{align}
\mathcal{P}:\:\min_{\boldsymbol{x},\Theta} & f_{0}(\boldsymbol{x},\Theta)\triangleq\mathbb{E}\left[g_{0}\left(\boldsymbol{x},\boldsymbol{y}\left(\boldsymbol{\xi}\right),\boldsymbol{\xi}\right)\right],\label{eq:mainP}\\
\text{s.t. } & f_{i}(\boldsymbol{x},\Theta)\triangleq\mathbb{E}\left[g_{i}\left(\boldsymbol{x},\boldsymbol{y}\left(\boldsymbol{\xi}\right),\boldsymbol{\xi}\right)\right]\leq0,\:i=1,...,m\nonumber \\
 & h_{j}\left(\boldsymbol{y}\left(\boldsymbol{\xi}\right),\boldsymbol{\xi}\right)\leq0,\:j=1,...,n,\forall\boldsymbol{\xi}\in\Omega,\label{eq:shortcon}
\end{align}
where $\boldsymbol{x}\in\mathcal{X}$ is the \textit{long-term (first
stage) optimization variable}, with $\mathcal{X}$ being the domain
of $\boldsymbol{x}$; $\boldsymbol{\xi}$ is a random state defined
on the probability space $\left(\Omega,\mathcal{F},\mathbb{P}\right)$,
with $\Omega$ being the sample space, $\mathcal{F}$ being the $\sigma$-algebra
generated by subsets of $\Omega$, and $\mathbb{P}$ being a probability
measure defined on $\mathcal{F}$; $\boldsymbol{y}\left(\boldsymbol{\xi}\right)$
is the \textit{short-term (second stage) optimization variable} under
state $\boldsymbol{\xi}$; and $\Theta\triangleq\left\{ \boldsymbol{y}\left(\boldsymbol{\xi}\right)\in\mathcal{Y},\forall\boldsymbol{\xi}\right\} $
is the collection of the short-term optimization variables for all
possible states, with $\mathcal{Y}$ being the domain of the short-term
optimization variable $\boldsymbol{y}$. $f_{i}(\boldsymbol{x},\Theta)\leq0,i=1,...,m$
are called the long-term constraints and $h_{j}\left(\boldsymbol{y}\left(\boldsymbol{\xi}\right),\boldsymbol{\xi}\right)\leq0,j=1,...,n$
are called the short-term constraints under the state $\boldsymbol{\xi}$.
Clearly, in $\mathcal{P}$, the long-term variable $\boldsymbol{x}$
is adaptive to the distribution/statistics of the random state $\boldsymbol{\xi}$,
while the short-term variable $\boldsymbol{y}$ is adaptive to the
realization of the random state $\boldsymbol{\xi}$.

The study of the stochastic optimization in Problem $\mathcal{P}$
is motivated by the following observations. Many physical systems
are not deterministic and we must take into account of the underlying
random state $\boldsymbol{\xi}$ in modelling optimization problems.
Moreover, in practice, some optimization variables have to be optimized
before observing the realization of the state $\boldsymbol{\xi}$,
while the other optimization variables can be optimized after observing
the realization of the state $\boldsymbol{\xi}$. These two types
of optimization variables can be modeled by the long-term optimization
variable and short-term optimization variable, respectively. For example,
consider a multi-access channel with $K$ users. There are $M$ antennas
at the receiver and a single antenna at each user (transmitter). The
receiver applies linear beamforming to decode the information, and
the average data rate of user $k$ is given by $\mathbb{E}\left[\log\left(1+\frac{p_{k}\left|\boldsymbol{u}_{k}^{H}\left(\boldsymbol{a}\right)\boldsymbol{a}_{k}\right|}{1+\sum_{l\neq k}p_{l}\left|\boldsymbol{u}_{k}^{H}\left(\boldsymbol{a}\right)\boldsymbol{a}_{l}\right|}\right)\right]$,
where $\boldsymbol{a}_{k}$ is the channel of user $k$, $\boldsymbol{a}=\left\{ \boldsymbol{a}_{k},\forall k\right\} $
is the aggregate channel state, $p_{k}$ is the transmit power for
user $k$, and $\boldsymbol{u}_{k}\left(\boldsymbol{a}\right)$ is
the receive beamforming vector of user $k$ for given channel state
$\boldsymbol{a}$. Suppose the users only know the channel statistics
and the receiver has perfect channel state information (CSI) $\boldsymbol{a}$.
As such, the transmit power $p_{k}$ is a long-term optimization variable
only adaptive to the channel statistics, and the receive beamforming
vectors $\boldsymbol{u}_{k}^{H}\left(\boldsymbol{a}\right),\forall\boldsymbol{a}$
are short optimization variables adaptive to the CSI $\boldsymbol{a}$.
The design goal is to minimize the total transmit power subject to
the average data rate constraints for each user:
\begin{align}
\min_{\left\{ p_{k},\boldsymbol{u}_{k}\left(\boldsymbol{a}\right)\right\} } & \sum_{k=1}^{K}p_{k}\label{eq:mainP-exm0}\\
\text{s.t. } & \mathbb{E}\left[\log\left(1+\frac{p_{k}\left|\boldsymbol{u}_{k}^{H}\left(\boldsymbol{a}\right)\boldsymbol{a}_{k}\right|}{1+\sum_{l\neq k}p_{l}\left|\boldsymbol{u}_{k}^{H}\left(\boldsymbol{a}\right)\boldsymbol{a}_{l}\right|}\right)\right]\leq0,\:\forall k.\nonumber 
\end{align}
Problem (\ref{eq:mainP-exm0}) is an instance of Problem $\mathcal{P}$
with random state $\boldsymbol{a}$. In fact, many important engineering
and scientific applications, such as those considered in wireless
resource optimizations \cite{Liu_TSP2014_BFRoutingD2D,Liu_TSP15_twostageCE,Liu_TSP18_CRAN},
transportation network design \cite{Hrabec2015} and machine/deep
learning \cite{Birge_Springer11_twostageSO,Lee_WCL2018_CRDNN,Lee_WCL2018_D2DDNN},
can be viewed as instances of the two-stage stochastic optimization
problem $\mathcal{P}$. Despite its wide applications, Problem $\mathcal{P}$
is very challenging and there only exist solutions for some specific
applications. The existing two-stage stochastic optimization algorithms
can be classified into the following four classes.

\textbf{Deterministic Algorithms based on Sample Average Approximation
(SAA):} In this class, problem $\mathcal{P}$ is approximated as a
deterministic optimization problem $\widetilde{\mathcal{P}}$ by replacing
the objective/constraint functions $f_{i}(\boldsymbol{x},\Theta),i=0,1,...,m$
with their SAAs $\widetilde{f}_{i}(\boldsymbol{x},\Theta)\triangleq\frac{1}{T}\sum_{j=1}^{T}g_{i}\left(\boldsymbol{x},\boldsymbol{y}\left(\boldsymbol{\xi}_{j}\right),\boldsymbol{\xi}_{j}\right),i=0,1,...,m$
using $T\gg1$ state samples $\boldsymbol{\xi}_{j},j=1,...,T$. Then,
various deterministic optimization algorithms such as majorization-minimization
(MM) \cite{Sun_TSP2017_MM} and successive convex approximation (SCA)
\cite{Scutari_TSP14_SCA} can be used to solve the resulting deterministic
optimization problem $\widetilde{\mathcal{P}}$ to obtain an approximate
solution for Problem $\mathcal{P}$. In order to achieve a good approximation,
$T$ is usually chosen to be a large number and thus the SAA-based
deterministic algorithms suffer from very high complexity.

\textbf{Primal-Dual Decomposition Algorithms:} The primal-dual methods
refer to the approaches which concurrently solving a primal problem
(corresponding to the original optimization task) as well as a dual
formulation of this problem \cite{Komodakis_SPM15_primaldual}. Primal-dual
methods have been primarily employed in convex optimization problems
and they are usually not guaranteed to converge in the non-convex
case \cite{Komodakis_SPM15_primaldual}. In \cite{Tanikawa_TAC85_Primaldual},
a nonconvex primal-dual decomposition method is proposed for \textit{separable
optimization problem}, where the optimization variables are separable
in both objective and constraint functions. However, it cannot be
applied to our problem because: 1) the long-term and short-term optimization
variables in $\mathcal{P}$ are coupled together in the objective/constraint
functions, which does not satisfy the separable assumption in \cite{Tanikawa_TAC85_Primaldual};
2) the nonconvex primal-dual decomposition method in \cite{Komodakis_SPM15_primaldual}
is only locally convergent to an (locally) optimal solution when the
initial point is sufficiently close to it. In \cite{Berkelaar_MP05_primaldual},
a new primal-dual decomposition algorithm is proposed for two-stage
stochastic optimization with a convex objective and stochastic recourse
matrices. However, it does not work for non-convex stochastic optimization.

\textbf{Two-stage Stochastic Algorithms with Increasing Batch Size:}
In this class, it is usually assumed that one independent state sample
can be observed at the beginning of each iteration of the algorithm.
Then, at the $t$-th iteration, the algorithm uses all the available
$t$ state samples observed in the previous $t$ iterations to update
the long-term variable and $t$ short-term variables associated with
the $t$ state samples, whose complexity may become unacceptable when
the number of iterations $t$ becomes large. A famous example of this
class is the stochastic cutting plane algorithm (SCPA), which only
works for two-stage stochastic convex problems \cite{higle1991stochastic,Liu_TSP2014_BFRoutingD2D}.
In \cite{Liu_TSP18_CRAN}, an approximate stochastic cutting plane
algorithm (ASCPA) is proposed to find a sub-optimal solution for a
class of two-stage stochastic non-convex optimization problems with
piece-wise linear objective functions. In \cite{Liu_TSP15_twostageCE},
an alternating optimization (AO) algorithm with increasing batch size
is proposed to find a stationary point for a two-stage stochastic
non-convex optimization problem.

\textbf{Two-stage Stochastic Algorithm with Finite Batch Size: }For
the special case when the long-term and short-term variables are decoupled
in the constraint, i.e., when $m=0$ and the coupled constraints $f_{i}(\boldsymbol{x},\Theta)\leq0$
is absent, a two-stage stochastic algorithm called TOSCA with finite
batch size is proposed in \cite{Liu_TSP2018_TOSCA} to find a stationary
point of Problem $\mathcal{P}$. In this case, it is assumed that
a finite batch of $B\geq1$ independent state samples can be observed
at the beginning of each iteration. Then, at each iteration, the TOSCA
algorithm uses $B$ state samples to update the long-term variable
and $B$ short-term variables associated with the $B$ state samples.
The batch size $B$ can be chosen to achieve a good tradeoff between
the per-iteration complexity and the convergence speed.

\textbf{Deep Learning (DL) Algorithms:} Recently, DL has been widely
applied to solve complicated optimization problems in various application
areas \cite{Hong_TSP2018_WMMSEDNN,Yu_JSAC2019_Deepscheduling,Debbah_TSP2020_DLpowcont}.
For example, \cite{Hong_TSP2018_WMMSEDNN} is one of the pioneer works
to consider resource allocation using DL-based optimization technique.
In \cite{Debbah_TSP2020_DLpowcont}, the implementability of the global
optimal solution with DL is demonstrated for a non-convex power allocation
problem.\textbf{ }A few DL algorithms have also been proposed to solve
$\mathcal{P}$ for the special case when the long-term variable $\boldsymbol{x}$
is absent and only the short-term variables $\boldsymbol{y}\left(\boldsymbol{\xi}\right)$
is present \cite{Lee_WCL2018_CRDNN,Lee_WCL2018_D2DDNN,Quek_JSAC2019_DLDist}.
In this case, we can use a DNN $\boldsymbol{\phi}\left(\boldsymbol{\xi};\boldsymbol{\theta}\right)$
to approximate the optimal solution $\boldsymbol{y}^{\star}\left(\boldsymbol{\xi}\right)$
and the optimal parameter $\boldsymbol{\theta}^{\star}$ for the DNN
can be found by solving the following unsupervised DL problem:
\begin{align}
\min_{\boldsymbol{\theta}}\: & \mathbb{E}\left[g_{0}\left(\boldsymbol{\phi}\left(\boldsymbol{\xi};\boldsymbol{\theta}\right),\boldsymbol{\xi}\right)\right],\label{eq:mainP-DNN}\\
\text{s.t. } & \mathbb{E}\left[g_{i}\left(\boldsymbol{\phi}\left(\boldsymbol{\xi};\boldsymbol{\theta}\right),\boldsymbol{\xi}\right)\right]\leq0,\:i=1,...,m\nonumber \\
 & h_{j}\left(\boldsymbol{\phi}\left(\boldsymbol{\xi};\boldsymbol{\theta}\right),\boldsymbol{\xi}\right)\leq0,\:j=1,...,n,\forall\boldsymbol{\xi}.\nonumber 
\end{align}
However, existing DL optimizers are not suitable for handling complicated
stochastic constraints in (\ref{eq:mainP-DNN}) \cite{Quek_JSAC2019_DLDist}.
To resolve this issue, a penalizing method is applied to transform
the original constrained training problem to an unconstrained one
by augmenting a penalty term associated with the complicated constraints
into the objective function \cite{Lee_WCL2018_CRDNN,Lee_WCL2018_D2DDNN}.
In \cite{Quek_JSAC2019_DLDist}, a primal-dual method is employed
to solve the constrained training problem in (\ref{eq:mainP-DNN}).
However, when both the objective and constraint functions in (\ref{eq:mainP-DNN})
are non-convex, the methods in \cite{Lee_WCL2018_CRDNN,Lee_WCL2018_D2DDNN,Quek_JSAC2019_DLDist}
cannot guarantee the convergence to a feasible stationary point of
the original problem $\mathcal{P}$. Moreover, the DNN solution $\boldsymbol{\phi}\left(\boldsymbol{\xi};\boldsymbol{\theta}\right)$
does not exploit the specific problem structure and the number of
parameters $\boldsymbol{\theta}$ is usually large, leading to a high
complexity and slow convergence of the DL algorithms.

In summary, there still lacks efficient algorithms for the two-stage
stochastic optimization problem $\mathcal{P}$, especially when the
long-term and short-term variables are tightly coupled in the constraints.
In this paper, we propose a two-stage stochastic optimization algorithmic
framework based on a primal-dual decomposition method and deep unrolling,
to overcome the disadvantages of the above existing algorithms. The
main contributions are summarized as follows.
\begin{itemize}
\item \textbf{Two-stage primal-dual decomposition method:} The existing
primal-dual (decomposition) methods cannot handle the tightly coupled
non-convex stochastic constraints containing both the long-term and
short-term variables, as explained above. To overcome this challenge,
we establish a novel two-stage primal-dual decomposition method to
decompose $\mathcal{P}$ into a family of short-term subproblems $\mathcal{P}_{S}\left(\boldsymbol{x},\boldsymbol{\lambda},\boldsymbol{\xi}\right),\forall\boldsymbol{\xi}$
for fixed long-term variable $\boldsymbol{x}$ and Lagrange multipliers
$\boldsymbol{\lambda}$, and a long-term problem $\mathcal{P}_{L}$
with both $\boldsymbol{x}$ and $\boldsymbol{\lambda}$ as the optimization
variables, where $\boldsymbol{\lambda}=\left[\lambda_{1},...,\lambda_{m}\right]^{T}\succeq\mathbf{0}$
are the Lagrange multipliers associated with the long-term constraints
$f_{i}(\boldsymbol{x},\Theta)\leq0,\:i=1,...,m$. Since this method
involves solving both a primal problem and a ``mixed-primal-dual''
problem (i.e., a problem containing both primal variable $\boldsymbol{x}$
and dual variable $\boldsymbol{\lambda}$), and the primal problem
is further decomposed into short-term subproblems, we call it the
two-stage ``primal-dual decomposition'' method. Such a decomposition
method is clearly different from the existing primal-dual methods
in \cite{Komodakis_SPM15_primaldual,Tanikawa_TAC85_Primaldual,Berkelaar_MP05_primaldual}.
We establish the global and local optimality of this decomposition
method in Theorem 1 and 2, respectively, whose proofs are non-trivial
as shown in Appendix \ref{subsec:Proof-of-Theorem-PDDidea} and \ref{subsec:Proof-of-Theorem-PDDprac}. 
\item \textbf{A two-stage stochastic optimization framework:} We propose
a primal-dual decomposition based stochastic successive convex approximation
(PDD-SSCA) framework which can be applied to a class of two-stage
stochastic optimization problems that satisfy certain smooth conditions
as will be given in Assumption 1, without any other restrictive assumptions.
For example, the long-term and short-term variables can have tight
coupling in the constraints and all objective and constraint functions
can be non-convex. As such, this optimization framework opens the
door to solving two-stage stochastic optimization problems that occur
in many new applications. We establish the convergence of PDD-SSCA
to KKT solutions of $\mathcal{P}$ under mild conditions. To the best
of our knowledge, the proposed PDD-CSSCA is the first algorithm that
can guarantee the convergence to a KKT solution of a two-stage stochastic
optimization problem with tightly coupled non-convex stochastic constraints. 
\item \textbf{Specific PDD-SSCA algorithm design for some important applications:}
We apply PDD-SSCA to solve two important problems in wireless resource
allocation and hybrid analog-digital signal processing, respectively.
We believe that the proposed PDD-SSCA solutions for these problems
alone are of great interest to the community.
\end{itemize}

The rest of the paper is organized as follows. The assumptions on
the problem formulation is given in Section \ref{sec:System-Model},
together with some application examples. The primal-dual decomposition
method for two-stage stochastic optimization is established in Section
\ref{sec:Theory-of-Primal-Dual}. The PDD-SSCA algorithm and the convergence
analysis are presented in Section \ref{sec:Constrained-Stochastic-Successiv},
and some implementation details are discussed in Section \ref{sec:Deep-Neural-Network}.
Section \ref{sec:Applications} applies PDD-SSCA to solve two important
application problems. Finally, the conclusion is given in Section
\ref{sec:Conclusion}.

\section{Problem Formulation\label{sec:System-Model}}

\subsection{Assumptions on Problem $\mathcal{P}$}

We make the following assumptions on Problem $\mathcal{P}$.

\begin{assumption}[Assumptions on Problem $\mathcal{P}$]\label{asm:convP}$\:$
\begin{enumerate}
\item $\mathcal{X}\subseteq\mathbb{R}^{n_{x}}$ and $\mathcal{Y}\subseteq\mathbb{R}^{n_{y}}$
for some positive integers $n_{x}$ and $n_{y}$. Moreover, $\mathcal{X},\mathcal{Y}$
are compact and convex. 
\item The functions $g_{i}\left(\boldsymbol{x},\boldsymbol{y},\boldsymbol{\xi}\right),i=0,...,m$
are real valued and continuously differentiable functions in $\boldsymbol{x}\in\mathcal{X},\boldsymbol{y}\in\mathcal{Y}$.
\item The functions $h_{j}\left(\boldsymbol{y},\boldsymbol{\xi}\right),j=0,...,n$
are real valued and continuously differentiable functions in $\boldsymbol{y}\in\mathcal{Y}$. 
\item For any $i\in\left\{ 0,...,m\right\} $ and $\boldsymbol{\xi}\in\Omega$,
the function $g_{i}\left(\boldsymbol{x},\boldsymbol{y},\boldsymbol{\xi}\right)$,
its derivative w.r.t. $\boldsymbol{x}$ and $\boldsymbol{y}$, and
its second-order derivative w.r.t. $\boldsymbol{x}$ and $\boldsymbol{y}$,
are uniformly bounded.
\item For any $j\in\left\{ 1,...,n\right\} $ and $\boldsymbol{\xi}\in\Omega$,
the function $h_{j}\left(\boldsymbol{y},\boldsymbol{\xi}\right)$,
its derivative w.r.t. $\boldsymbol{y}$, and its second-order derivative
w.r.t. $\boldsymbol{y}$, are uniformly bounded.
\item For any $\boldsymbol{\xi}\in\Omega$, the short-term constraints are
feasible, i.e., $\exists\boldsymbol{y}\in\mathcal{Y}$, such that
$h_{j}\left(\boldsymbol{y},\boldsymbol{\xi}\right)\leq0,\forall j$.
Moreover, Problem $\mathcal{P}$ is feasible.
\end{enumerate}
\end{assumption}

These conditions are standard and are satisfied for a large class
of problems. For ease of exposition, we assume $\boldsymbol{x},\boldsymbol{y}$
are real vectors. Nevertheless, the proposed algorithm can be directly
applied to the case with complex optimization variables $\boldsymbol{x},\boldsymbol{y}$,
by treating each function $g_{i}\left(\boldsymbol{x},\boldsymbol{y},\xi\right)$
in the problem as a real valued function of real vectors $\left[\textrm{Re}\left[\boldsymbol{x}\right];\textrm{Im}\left[\boldsymbol{x}\right]\right]$
and $\left[\textrm{Re}\left[\boldsymbol{y}\right];\textrm{Im}\left[\boldsymbol{y}\right]\right]$.

\subsection{Examples of Problem $\mathcal{P}$\label{subsec:Examples-of-Problem}}

Problem $\mathcal{P}$ embraces many applications. In the following,
we give two important examples. 
\begin{example}
[Cognitive Multiple Access Channels \cite{Zhang_TIT2009_CMAC}]\label{exa:Massive-MIMO-Hybrid}
Consider a multi-user uplink cognitive (CR) network \cite{Zhang_TIT2009_CMAC}
with one licensed primary user (PU) and $N$ secondary users (SUs).
The SUs share time and frequency resources with the PU and desire
to transmit their data to a secondary base station (SBS). The SBS
and all users are all equipped with only a single antenna. Let $a_{i}$
and $b_{i}$ denote the channel gain from SU $i$ to the SBS and the
PU, respectively. The average sum capacity maximization problem can
be formulated as \cite{Zhang_TIT2009_CMAC}
\begin{align}
\max_{\left\{ \boldsymbol{p}\left(\boldsymbol{a},\boldsymbol{b}\right)\right\} } & \mathbb{E}\left[\log\left(1+\sum_{i=1}^{N}a_{i}p_{i}\left(\boldsymbol{a},\boldsymbol{b}\right)\right)\right],\label{eq:mainP-exm1}\\
\text{s.t. } & \mathbb{E}\left[p_{i}\left(\boldsymbol{a},\boldsymbol{b}\right)\right]\leq P_{i},\:i=1,...,N\nonumber \\
 & \mathbb{E}\left[\sum_{i=1}^{N}b_{i}p_{i}\left(\boldsymbol{a},\boldsymbol{b}\right)\right]\leq\Gamma,\nonumber \\
 & p_{i}\left(\boldsymbol{a},\boldsymbol{b}\right)\geq0,\:i=1,...,N,\forall\boldsymbol{a},\boldsymbol{b},\nonumber 
\end{align}
where $\boldsymbol{a}=\left[a_{1},...,a_{N}\right]^{T}$, $\boldsymbol{b}=\left[b_{1},...,b_{N}\right]^{T}$,
$p_{i}\left(\boldsymbol{a},\boldsymbol{b}\right)$ is the transmit
power at SU $i$ when the channel state is $\boldsymbol{a},\boldsymbol{b}$,
$\boldsymbol{p}\left(\boldsymbol{a},\boldsymbol{b}\right)=\left[p_{1}\left(\boldsymbol{a},\boldsymbol{b}\right),...,p_{N}\left(\boldsymbol{a},\boldsymbol{b}\right)\right]^{T}$,
$P_{i}$ and $\Gamma$ stand for the long-term transmit power budget
at SU $i$ and the interference threshold constraint for the PU, respectively.
Clearly, problem (\ref{eq:mainP-exm1}) is a special case of $\mathcal{P}$
with short-term variables $\boldsymbol{p}\left(\boldsymbol{a},\boldsymbol{b}\right),\forall\boldsymbol{a},\boldsymbol{b}$
only. The random state is the channel state $\boldsymbol{a},\boldsymbol{b}$.
\end{example}

\begin{example}
[Power Minimization for Two-timescale hybrid beamforming \cite{Liu_TSP14_RFprecoding}]\label{exa:Massive-MIMO-Hybrid-1}Massive
MIMO is considered as one of the key technologies in 5G wireless systems.
Consider a multi-user massive MIMO downlink system where a base station
(BS) serves $K$ single-antenna users. The BS is equipped with $M\gg1$
antennas and $S$ transmit RF chains, where $K\leq S<M$. Two-timescale
hybrid beamforming is employed at the BS to support simultaneous transmissions
to the $K$ users, with reduced hardware cost and channel state information
(CSI) signaling overhead \cite{Liu_TSP14_RFprecoding,Liu_TSP2016_CSImassive,Park_TSP17_THP}.
Specifically, the precoder is split into a baseband precoder and an
RF precoder as $\boldsymbol{F}\boldsymbol{G}$, where $\boldsymbol{G}=\left[\boldsymbol{g}_{1},...,\boldsymbol{g}_{K}\right]\in\mathbb{C}^{S\times K}$
is the baseband precoder, and $\boldsymbol{F}\in\mathbb{C}^{M\times S}$
is the RF precoder using the RF phase shifting network \cite{Zhang_TSP05_RFshifter}.
Hence, all elements of $\boldsymbol{F}$ have equal magnitude, i.e.,
$F_{m,s}=e^{j\theta_{m,s}}$, where $\theta_{m,s}$ is the phase of
the $\left(m,s\right)$-th element $F_{m,s}$ of $\boldsymbol{F}$.
The RF precoder $\boldsymbol{F}$ can be represented by a phase vector
$\boldsymbol{\theta}\in\mathbb{R}^{MS}$ whose $\left(\left(j-1\right)M+i\right)$-th
element is $\theta_{i,j}$.

We focus on a coherence time interval of channel statistics within
which the channel statistics (distribution) are assumed to be constant.
The coherence time of channel statistics is divided into $T_{f}$
frames and each frame consists of $T_{s}$ time slots, as illustrated
in Fig. \ref{fig:timelineAlg}. The channel state $\boldsymbol{H}=\left[\boldsymbol{h}_{1},...,\boldsymbol{h}_{K}\right]^{H}\in\mathbb{C}^{K\times M}$
is assumed to be constant within each time slot, where $\boldsymbol{h}_{k}\in\mathbb{C}^{M}$
is the channel vector of user $k$. We assume that the BS can obtain
the real-time effective CSI $\widetilde{\boldsymbol{H}}=\boldsymbol{H}\boldsymbol{F}\in\mathbb{C}^{K\times S}$
at each time slot, and one outdated channel sample $\boldsymbol{H}$
at each frame. In the THP design, the analog precoder $\mathbf{F}$
is only updated once per frame based on the outdated channel sample
$\boldsymbol{H}$ to achieve massive MIMO array gain. The digital
precoder $\boldsymbol{G}$ is adaptive to the real-time effective
CSI $\widetilde{\boldsymbol{H}}\in\mathbb{C}^{K\times S}$ to achieve
the spatial multiplexing gain. Note that the effective CSI $\widetilde{\boldsymbol{H}}$
usually has much lower dimension than the full channel sample $\boldsymbol{H}$.
Therefore, it is possible to obtain the real-time effective CSI $\widetilde{\boldsymbol{H}}$
at each time slot by sending pilot signals with analog precoder $\boldsymbol{F}$.
However, we can only obtain one outdated full channel sample $\boldsymbol{H}$
at each frame because obtaining the real-time full CSI $\boldsymbol{H}$
per time slot will cause unacceptable CSI signaling overhead in massive
MIMO. Therefore, we cannot optimize both the analog and digital precoders
based on the real-time full CSI $\boldsymbol{H}$ at each time slot.
$\boldsymbol{F}$ and $\boldsymbol{G}$ have to be optimized at different
timescale based on the outdated channel sample $\boldsymbol{H}$ and
real-time effective CSI $\widetilde{\boldsymbol{H}}$, respectively,
using e.g., the proposed PDD-SSCA algorithm. In this example, each
frame corresponds to an iteration of PDD-SSCA.

\begin{figure}
\begin{centering}
\includegraphics[width=88mm]{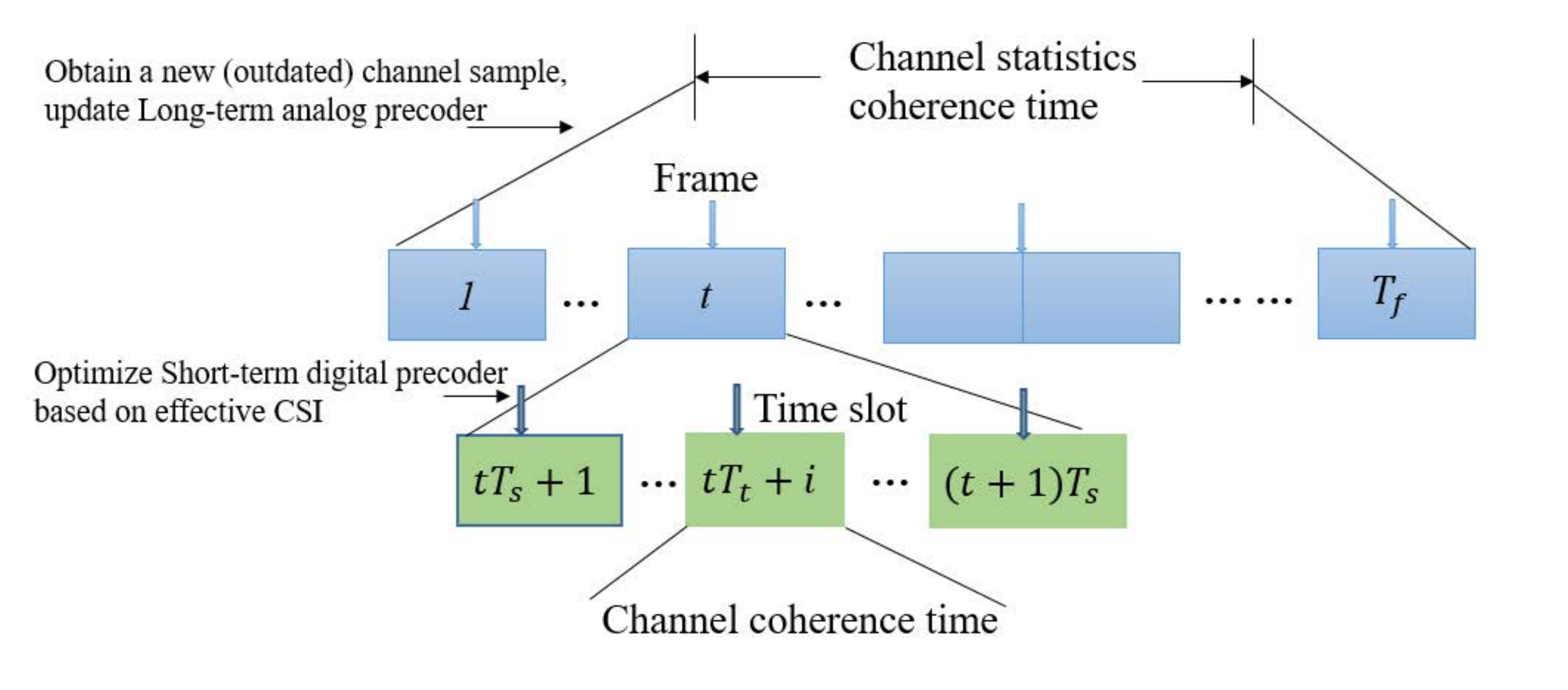}
\par\end{centering}
\caption{\label{fig:timelineAlg}Timeline (frame structure) of two-timescale
hybrid beamforming}
\end{figure}

For given RF precoding phase vector $\boldsymbol{\theta}$, baseband
precoder $\boldsymbol{G}$ and channel realization $\boldsymbol{H}$,
the data rate of user $k$ is given by
\begin{equation}
r_{k}\left(\boldsymbol{\theta},\boldsymbol{G},\boldsymbol{H}\right)=\log\left(1+\frac{\left|\boldsymbol{h}_{k}^{H}\boldsymbol{F}\boldsymbol{g}_{k}\right|^{2}}{\sum_{i\neq k}\left|\boldsymbol{h}_{k}^{H}\boldsymbol{F}\boldsymbol{g}_{i}\right|^{2}+1}\right).\label{eq:ratek}
\end{equation}
Note that $\boldsymbol{F}$ is a function of $\boldsymbol{\theta}$.
Consider the problem of average transmit power minimization for the
above massive MIMO system with individual average rate constraint
for each user, which can be formulated as:

\begin{align}
\min_{\boldsymbol{\theta},\left\{ \boldsymbol{G}\left(\boldsymbol{H}\right),\forall\boldsymbol{H}\right\} } & \mathbb{E}\left[Tr\left(\boldsymbol{F}\boldsymbol{G}\left(\boldsymbol{H}\right)\boldsymbol{G}^{H}\left(\boldsymbol{H}\right)\boldsymbol{F}^{H}\right)\right],\label{eq:NSP}\\
\text{s.t. } & \mathbb{E}\left[r_{k}\left(\boldsymbol{\theta},\boldsymbol{G}\left(\boldsymbol{H}\right),\boldsymbol{H}\right)\right]\geq\gamma_{k},k=1,...,K,\nonumber 
\end{align}
where $\gamma_{k}$ is the throughput requirement for user $k$. Problem
(\ref{eq:NSP}) is an instance of Problem $\mathcal{P}$ with random
state $\boldsymbol{H}$.
\end{example}

\subsection{KKT Solution of Problem $\mathcal{P}$}

Since Problem $\mathcal{P}$ is in general non-convex, we focus on
designing an efficient algorithm to find KKT solutions of Problem
$\mathcal{P}$, defined as follows.
\begin{defn}
[KKT solution of $\mathcal{P}$]\label{def:Stationary-solution}A
solution $\left(\boldsymbol{x}^{*}\in\mathcal{X},\Theta^{*}=\left\{ \boldsymbol{y}^{*}\left(\boldsymbol{\xi}\right)\in\mathcal{Y},\forall\boldsymbol{\xi}\right\} \right)$
is called a KKT solution of Problem $\mathcal{P}$, if there exist
long-term Lagrange multipliers $\boldsymbol{\lambda}=\left[\lambda_{1},...,\lambda_{m}\right]^{T}\succeq\mathbf{0}$
associated with the long-term constraints and short-term Lagrange
multipliers $\nu_{j}\left(\boldsymbol{\xi}\right)\geq0,\forall j,\forall\boldsymbol{\xi}$
associated with the short-term constraints, such that the following
conditions are satisfied: 
\begin{enumerate}
\item For every $\boldsymbol{\xi}\in\Omega$ outside a set of probability
zero, we have
\begin{align}
 & \partial_{\boldsymbol{y}}g_{0}\left(\boldsymbol{x}^{*},\boldsymbol{y}^{*}\left(\boldsymbol{\xi}\right),\boldsymbol{\xi}\right)\nonumber \\
 & +\sum_{i}\lambda_{i}\partial_{\boldsymbol{y}}g_{i}\left(\boldsymbol{x}^{*},\boldsymbol{y}^{*}\left(\boldsymbol{\xi}\right),\boldsymbol{\xi}\right)\nonumber \\
 & +\sum_{j}\nu_{j}\left(\boldsymbol{\xi}\right)\partial_{\boldsymbol{y}}h_{j}\left(\boldsymbol{y}^{*}\left(\boldsymbol{\xi}\right),\boldsymbol{\xi}\right)=\boldsymbol{0},\nonumber \\
 & h_{j}\left(\boldsymbol{y}^{*}\left(\boldsymbol{\xi}\right),\boldsymbol{\xi}\right)\leq0,\:j=1,...,n,\nonumber \\
 & \nu_{j}\left(\boldsymbol{\xi}\right)h_{j}\left(\boldsymbol{y}^{*}\left(\boldsymbol{\xi}\right),\boldsymbol{\xi}\right)=0,\:j=1,...,n.\label{eq:STscon}
\end{align}
where $\partial_{\boldsymbol{y}}g_{i}\left(\boldsymbol{x}^{*},\boldsymbol{y}^{*}\left(\boldsymbol{\xi}\right),\boldsymbol{\xi}\right)$
and $\partial_{\boldsymbol{y}}h_{j}\left(\boldsymbol{y}^{*}\left(\boldsymbol{\xi}\right),\boldsymbol{\xi}\right)$
are the partial derivatives of $g_{i}\left(\boldsymbol{x}^{*},\boldsymbol{y},\boldsymbol{\xi}\right)$
and $h_{j}\left(\boldsymbol{y},\boldsymbol{\xi}\right)$ w.r.t. $\boldsymbol{y}$
at $\boldsymbol{y}=\boldsymbol{y}^{*}\left(\boldsymbol{\xi}\right)$,
respectively.
\item 
\begin{align}
\partial_{\boldsymbol{x}}f_{0}\left(\boldsymbol{x}^{*},\Theta^{*}\right)+\sum_{i}\lambda_{i}\partial_{\boldsymbol{x}}f_{i}(\boldsymbol{x}^{*},\Theta^{*}) & =\boldsymbol{0},\nonumber \\
f_{i}(\boldsymbol{x}^{*},\Theta^{*}) & \leq0,\:\forall i\label{eq:LTscon}
\end{align}
where $\partial_{\boldsymbol{x}}f_{i}(\boldsymbol{x}^{*},\Theta^{*})$
is the partial derivative of $f_{i}(\boldsymbol{x},\Theta^{*})$ w.r.t.
$\boldsymbol{x}$ at $\boldsymbol{x}=\boldsymbol{x}^{*}$.
\item 
\begin{align}
 & \lambda_{i}f_{i}(\boldsymbol{x}^{*},\Theta^{*})=0,\:i=1,...,m.\label{eq:CompS}
\end{align}
\end{enumerate}
\end{defn}

\section{Two-Stage Primal-Dual Decomposition Method\label{sec:Theory-of-Primal-Dual}}

\subsection{Two-Stage Primal-Dual Decomposition for $\mathcal{P}$ }

One major challenge of solving $\mathcal{P}$ is that, the long-term
variable $\boldsymbol{x}$ and the short-term variables $\boldsymbol{y}\left(\boldsymbol{\xi}\right)$'s
for different states are coupled together in a complicated manner
via the long-term constraint. As discussed in the introduction, the
existing primal-dual (decomposition) methods in \cite{Komodakis_SPM15_primaldual,Tanikawa_TAC85_Primaldual,Berkelaar_MP05_primaldual}
cannot work for problem $\mathcal{P}$ with tightly coupled non-convex
stochastic constraints. To overcome this challenge, we prove a novel
primal-dual decomposition method to decouple the optimization variables.
Specifically, for a fixed long-term variable $\boldsymbol{x}$ and
long-term Lagrange multipliers $\boldsymbol{\lambda}$, let $\phi_{\boldsymbol{x},\boldsymbol{\lambda}}^{\star}$
denote the \textit{optimal short-term policy}, which is defined a
mapping from $\Omega$ to $\mathcal{Y}$ such that $\boldsymbol{y}^{\star}\left(\boldsymbol{x},\boldsymbol{\lambda},\boldsymbol{\xi}\right)=\phi_{\boldsymbol{x},\boldsymbol{\lambda}}^{\star}\left(\boldsymbol{\xi}\right)$
is the optimal solution of the following \textit{short-term subproblem}:
\begin{align}
\mathcal{P}_{S}\left(\boldsymbol{x},\boldsymbol{\lambda},\boldsymbol{\xi}\right):\:\min_{\boldsymbol{y}} & g_{0}\left(\boldsymbol{x},\boldsymbol{y},\boldsymbol{\xi}\right)+\sum_{i}\lambda_{i}g_{i}\left(\boldsymbol{x},\boldsymbol{y},\boldsymbol{\xi}\right),\label{eq:mainPS}\\
\text{s.t. } & h_{j}\left(\boldsymbol{y},\boldsymbol{\xi}\right)\leq0,\:j=1,...,n.\nonumber 
\end{align}
Note that the optimal short-term policy is not necessarily unique
and the set of all optimal short-term policies is denoted as $\Phi_{\boldsymbol{x},\boldsymbol{\lambda}}^{\star}$.
When $\Phi_{\boldsymbol{x},\boldsymbol{\lambda}}^{\star}$ have multiple
elements, we choose one optimal short-term policy $\phi_{\boldsymbol{x},\boldsymbol{\lambda}}^{\star}$
as $\phi_{\boldsymbol{x},\boldsymbol{\lambda}}^{\star}\in\text{argmin}_{\phi\in\Phi_{\boldsymbol{x},\boldsymbol{\lambda}}^{F}}\max_{i\in\left\{ 1,...,m\right\} }\mathbb{E}\left[g_{i}\left(\boldsymbol{x},\phi\left(\boldsymbol{\xi}\right),\boldsymbol{\xi}\right)\right]$,
where $\Phi_{\boldsymbol{x},\boldsymbol{\lambda}}^{F}\triangleq\text{argmin}_{\phi\in\Phi_{\boldsymbol{x},\boldsymbol{\lambda}}^{\star}}\sum_{i}\lambda_{i}\left|\mathbb{E}\left[g_{i}\left(\boldsymbol{x},\phi\left(\boldsymbol{\xi}\right),\boldsymbol{\xi}\right)\right]\right|$
is the set of all optimal short-term policies that minimize $\sum_{i}\lambda_{i}\left|\mathbb{E}\left[g_{i}\left(\boldsymbol{x},\phi\left(\boldsymbol{\xi}\right),\boldsymbol{\xi}\right)\right]\right|$.
In other words, we choose $\phi_{\boldsymbol{x},\boldsymbol{\lambda}}^{\star}$
as an optimal short-term policy which is most likely to satisfy the
constraints and complementary slackness condition of the original
Problem $\mathcal{P}$. With the optimal short-term policy $\phi_{\boldsymbol{x},\boldsymbol{\lambda}}^{\star}$
and $\left\{ \boldsymbol{y}^{\star}\left(\boldsymbol{x},\boldsymbol{\lambda},\boldsymbol{\xi}\right)=\phi_{\boldsymbol{x},\boldsymbol{\lambda}}^{\star}\left(\boldsymbol{\xi}\right),\forall\boldsymbol{\xi}\right\} $
chosen according to the above rule, we formulate the following long-term
problem
\begin{align}
\mathcal{P}_{L}:\:\min_{\boldsymbol{x},\boldsymbol{\lambda}} & f_{0}^{\star}(\boldsymbol{x},\boldsymbol{\lambda})\triangleq\mathbb{E}\left[g_{0}\left(\boldsymbol{x},\boldsymbol{y}^{\star}\left(\boldsymbol{x},\boldsymbol{\lambda},\boldsymbol{\xi}\right),\boldsymbol{\xi}\right)\right],\label{eq:mainP-1-1}\\
\text{s.t. } & f_{i}^{\star}(\boldsymbol{x},\boldsymbol{\lambda})\triangleq\mathbb{E}\left[g_{i}\left(\boldsymbol{x},\boldsymbol{y}^{\star}\left(\boldsymbol{x},\boldsymbol{\lambda},\boldsymbol{\xi}\right),\boldsymbol{\xi}\right)\right]\leq0,\:\forall i.\nonumber 
\end{align}
Problem $\mathcal{P}_{L}$ only contains long-term variables $\boldsymbol{x},\boldsymbol{\lambda}$.
Now we are ready to establish a two-stage primal-dual decomposition
theorem for Problem $\mathcal{P}$.
\begin{thm}
[Two-stage Primal-Dual Decomposition]\label{thm:PDidea}Let $\left(\boldsymbol{x}^{\star},\boldsymbol{\lambda}^{\star}\right)$
denote any optimal solution of $\mathcal{P}_{L}$. Define 
\begin{align*}
G_{i}^{\textrm{min}}\left(\boldsymbol{x}\right) & \triangleq\min_{\Theta}f_{i}(\boldsymbol{x},\Theta),\text{ s.t. (\ref{eq:shortcon}) is satisfied},\\
G_{i}^{\textrm{max}}\left(\boldsymbol{x}\right) & \triangleq\max_{\Theta}f_{i}(\boldsymbol{x},\Theta)\text{ s.t. (\ref{eq:shortcon}) is satisfied}.
\end{align*}
If there exists an optimal solution of $\mathcal{P}$, denoted as
$\boldsymbol{x}^{\circ},\Theta^{\circ}=\left\{ \boldsymbol{y}^{\circ}\left(\boldsymbol{\xi}\right),\forall\boldsymbol{\xi}\right\} $,
such that there exist arbitrary small numbers $\delta_{i}\in\left(0,G_{i}^{\textrm{max}}\left(\boldsymbol{x}^{\circ}\right)-G_{i}^{\textrm{min}}\left(\boldsymbol{x}^{\circ}\right)\right),i=1,...,m$,
then $\left(\boldsymbol{x}^{\star},\Theta^{\star}=\left\{ \boldsymbol{y}^{\star}\left(\boldsymbol{x}^{\star},\boldsymbol{\lambda}^{\star},\boldsymbol{\xi}\right)\in\mathcal{Y},\forall\boldsymbol{\xi}\right\} \right)$
is also the optimal solution of $\mathcal{P}$.
\end{thm}

Please refer to Appendix \ref{subsec:Proof-of-Theorem-PDDidea} for
the proof.

Note that the condition in Theorem \ref{thm:PDidea} means that for
fixed long-term variable $\boldsymbol{x}^{\circ}$, the maximum value
of the constraint function $G_{i}^{\textrm{max}}\left(\boldsymbol{x}^{\circ}\right)$
is strictly larger than the minimum value of the constraint function
$G_{i}^{\textrm{min}}\left(\boldsymbol{x}^{\circ}\right)$, which
can be easily satisfied in practice. Theorem \ref{thm:PDidea} essentially
states that $\mathcal{P}$ can be decomposed into a family of short-term
subproblems $\mathcal{P}_{S}\left(\boldsymbol{x},\boldsymbol{\lambda},\boldsymbol{\xi}\right)$
for fixed $\boldsymbol{x},\boldsymbol{\lambda},\boldsymbol{\xi}$
and a long-term problem $\mathcal{P}_{L}$ with $\boldsymbol{x},\boldsymbol{\lambda}$
as optimization variables. However, it is not convenient to directly
apply Theorem \ref{thm:PDidea} for algorithm design because such
a decomposition requires the optimal solution of the short-term subproblem
$\mathcal{P}_{S}\left(\boldsymbol{x},\boldsymbol{\lambda},\boldsymbol{\xi}\right)$
for each $\boldsymbol{\xi}$, which is difficult to obtain in practice,
especially when $\mathcal{P}_{S}\left(\boldsymbol{x},\boldsymbol{\lambda},\boldsymbol{\xi}\right)$
is non-convex. Therefore, in the following, we will establish a relaxed
primal-dual decomposition method which does not require the optimal
solution of $\mathcal{P}_{S}\left(\boldsymbol{x},\boldsymbol{\lambda},\boldsymbol{\xi}\right)$.
Consequently, the relaxed primal-dual decomposition method can be
applied to design an efficient algorithm to find KKT solutions of
$\mathcal{P}$ up to certain tolerable error.

\subsection{Relaxed Two-Stage Primal-Dual Decomposition }

Let $\boldsymbol{y}^{J}\left(\boldsymbol{x},\boldsymbol{\lambda},\boldsymbol{\xi}\right)$
denote a stationary point (up to certain tolerable error) of the \textit{short-term
subproblem} $\mathcal{P}_{S}\left(\boldsymbol{x},\boldsymbol{\lambda},\boldsymbol{\xi}\right)$
obtained by running a \textit{short-term sub-algorithm} for $J$ iterations.
The short-term sub-algorithm is basically an iterative algorithm to
find a stationary point of the short-term subproblem $\mathcal{P}_{S}\left(\boldsymbol{x},\boldsymbol{\lambda},\boldsymbol{\xi}\right)$.
In other words, $\boldsymbol{y}^{J}\left(\boldsymbol{x},\boldsymbol{\lambda},\boldsymbol{\xi}\right)$
satisfies the KKT conditions of $\mathcal{P}_{S}\left(\boldsymbol{x},\boldsymbol{\lambda},\boldsymbol{\xi}\right)$
(up to certain tolerable error) as
\begin{align}
 & \left\Vert \partial_{\boldsymbol{y}}g_{0}\left(\boldsymbol{x},\boldsymbol{y}^{J}\left(\boldsymbol{x},\boldsymbol{\lambda},\boldsymbol{\xi}\right),\boldsymbol{\xi}\right)\right.\nonumber \\
 & +\sum_{i}\lambda_{i}\partial_{\boldsymbol{y}}g_{i}\left(\boldsymbol{x},\boldsymbol{y}^{J}\left(\boldsymbol{x},\boldsymbol{\lambda},\boldsymbol{\xi}\right),\boldsymbol{\xi}\right)\nonumber \\
 & +\left.\sum_{j}\nu_{j}\left(\boldsymbol{x},\boldsymbol{\lambda},\boldsymbol{\xi}\right)\partial_{\boldsymbol{y}}h_{j}\left(\boldsymbol{y}^{J}\left(\boldsymbol{x},\boldsymbol{\lambda},\boldsymbol{\xi}\right),\boldsymbol{\xi}\right)\right\Vert =e_{1}^{J}\left(\boldsymbol{x},\boldsymbol{\lambda},\boldsymbol{\xi}\right),\nonumber \\
 & h_{j}\left(\boldsymbol{y}^{J}\left(\boldsymbol{x},\boldsymbol{\lambda},\boldsymbol{\xi}\right),\boldsymbol{\xi}\right)\leq e_{2}^{J}\left(\boldsymbol{x},\boldsymbol{\lambda},\boldsymbol{\xi}\right),\:j=1,...,n,\nonumber \\
 & \nu_{j}(\boldsymbol{x},\boldsymbol{\lambda},\boldsymbol{\xi})h_{j}\left(\boldsymbol{y}^{J}(\boldsymbol{x},\boldsymbol{\lambda},\boldsymbol{\xi}),\boldsymbol{\xi}\right)=e_{3,j}^{J}(\boldsymbol{x},\boldsymbol{\lambda},\boldsymbol{\xi}),\:j=1,...,n,\label{eq:STsconeJ}
\end{align}
where $\nu_{j}\left(\boldsymbol{x},\boldsymbol{\lambda},\boldsymbol{\xi}\right)\geq0,\forall j$
are the short-term Lagrange multipliers, and $e_{1}^{J}\left(\boldsymbol{x},\boldsymbol{\lambda},\boldsymbol{\xi}\right),e_{2}^{J}\left(\boldsymbol{x},\boldsymbol{\lambda},\boldsymbol{\xi}\right),e_{3,j}^{J}(\boldsymbol{x},\boldsymbol{\lambda},\boldsymbol{\xi}),j=1,...,n$
is the error due to that the short-term sub-algorithm only runs for
a finite number of $J$ iterations. Suppose that the short-term sub-algorithm
converges to a stationary point of $\mathcal{P}_{S}\left(\boldsymbol{x},\boldsymbol{\lambda},\boldsymbol{\xi}\right)$
as $J\rightarrow\infty$. Then for all $\boldsymbol{x}\in\mathcal{X},\boldsymbol{\lambda}\succeq\mathbf{0}$,
we have $\lim_{J\rightarrow\infty}e_{i}^{J}\left(\boldsymbol{x},\boldsymbol{\lambda},\boldsymbol{\xi}\right)=0,i=1,2$,
$\lim_{J\rightarrow\infty}e_{3,j}^{J}\left(\boldsymbol{x},\boldsymbol{\lambda},\boldsymbol{\xi}\right)=0,\forall j$
. Note that in practice, the short-term sub-algorithm always runs
for a finite number of iterations. Therefore, it is meaningful to
derive a relaxed primal-dual decomposition theorem for this case.

With the notation of $\boldsymbol{y}^{J}\left(\boldsymbol{x},\boldsymbol{\lambda},\boldsymbol{\xi}\right)$,
we formulate the following relaxed long-term problem
\begin{align}
\mathcal{P}_{L}^{J}:\:\min_{\boldsymbol{x},\boldsymbol{\lambda}} & f_{0}^{J}(\boldsymbol{x},\boldsymbol{\lambda})\triangleq\mathbb{E}\left[g_{0}\left(\boldsymbol{x},\boldsymbol{y}^{J}\left(\boldsymbol{x},\boldsymbol{\lambda},\boldsymbol{\xi}\right),\boldsymbol{\xi}\right)\right],\label{eq:mainP-1}\\
\text{s.t. } & f_{i}^{J}(\boldsymbol{x},\boldsymbol{\lambda})\triangleq\mathbb{E}\left[g_{i}\left(\boldsymbol{x},\boldsymbol{y}^{J}\left(\boldsymbol{x},\boldsymbol{\lambda},\boldsymbol{\xi}\right),\boldsymbol{\xi}\right)\right]\leq0,\:\forall i.\nonumber 
\end{align}
Problem $\mathcal{P}_{L}^{J}$ only contains long-term variables $\boldsymbol{x},\boldsymbol{\lambda}$.
Let $\nabla_{\boldsymbol{x}}f_{i}^{J}(\boldsymbol{x},\boldsymbol{\lambda})\triangleq\mathbb{E}\left[\partial_{\boldsymbol{x}}\boldsymbol{y}^{J}\partial_{\boldsymbol{y}}g_{i}\left(\boldsymbol{x},\boldsymbol{y}^{J},\boldsymbol{\xi}\right)+\partial_{\boldsymbol{x}}g_{i}\left(\boldsymbol{x},\boldsymbol{y}^{J},\boldsymbol{\xi}\right)\right]$
and $\nabla_{\boldsymbol{\lambda}}f_{i}^{J}(\boldsymbol{x},\boldsymbol{\lambda})\triangleq\mathbb{E}\left[\partial_{\boldsymbol{\lambda}}\boldsymbol{y}^{J}\partial_{\boldsymbol{y}}g_{i}\left(\boldsymbol{x},\boldsymbol{y}^{J},\boldsymbol{\xi}\right)\right]$
denote the derivative of $f_{i}^{J}(\boldsymbol{x},\boldsymbol{\lambda})$
w.r.t. $\boldsymbol{x}$ and $\boldsymbol{\lambda}$, respectively,
where $\partial_{\boldsymbol{x}}\boldsymbol{y}^{J}=\partial_{\boldsymbol{x}}\boldsymbol{y}^{J}\left(\boldsymbol{x},\boldsymbol{\lambda},\boldsymbol{\xi}\right)\in\mathbb{R}^{n_{x}\times n_{y}}$
($\partial_{\boldsymbol{\lambda}}\boldsymbol{y}^{J}=\partial_{\boldsymbol{\lambda}}\boldsymbol{y}^{J}\left(\boldsymbol{x},\boldsymbol{\lambda},\boldsymbol{\xi}\right)\in\mathbb{R}^{m\times n_{y}}$)
is the derivative of the vector function $\boldsymbol{y}^{J}\left(\boldsymbol{x},\boldsymbol{\lambda},\boldsymbol{\xi}\right)$
to the vector $\boldsymbol{x}$ ($\boldsymbol{\lambda}$), $\boldsymbol{y}^{J}$
is an abbreviation for $\boldsymbol{y}^{J}\left(\boldsymbol{x},\boldsymbol{\lambda},\boldsymbol{\xi}\right)$.
Note that throughout this paper, we use $\nabla_{\boldsymbol{x}},\nabla_{\boldsymbol{\lambda}}$
to denote the derivative of $f_{i}^{J}$ or $g_{i}$ by viewing $\boldsymbol{y}^{J}\left(\boldsymbol{x},\boldsymbol{\lambda},\boldsymbol{\xi}\right)$
as a function of $\boldsymbol{x},\boldsymbol{\lambda}$, and use $\partial_{\boldsymbol{x}},\partial_{\boldsymbol{\lambda}}$
to denote the derivative of $f_{i}^{J}$ or $g_{i}$ by viewing $\boldsymbol{y}^{J}\left(\boldsymbol{x},\boldsymbol{\lambda},\boldsymbol{\xi}\right)$
as a fixed value.

In the following theorem, we establish a relaxed primal-dual decomposition
theorem for the two-stage Problem $\mathcal{P}$, which provides a
foundation for the algorithm design.
\begin{thm}
[Relaxed Two-stage Primal-Dual Decomposition]\label{thm:PDprac}Suppose
that for every $\boldsymbol{\xi}\in\Omega$ outside a set of probability
zero, $\boldsymbol{y}^{J}\left(\boldsymbol{x},\boldsymbol{\lambda},\boldsymbol{\xi}\right)$
is continuously differentiable function in $\boldsymbol{x}\in\mathcal{X}$
and $\boldsymbol{\lambda}\succeq\mathbf{0}$. Let $\boldsymbol{x}^{*},\boldsymbol{\lambda}^{*}$
denote a KKT point of $\mathcal{P}_{L}^{J}$, i.e., there exists Lagrange
multipliers $\widetilde{\boldsymbol{\lambda}}=\left[\widetilde{\lambda}_{1},...,\widetilde{\lambda}_{m}\right]^{T}\succeq\mathbf{0}$
such that the following KKT conditions are satisfied:
\begin{align}
\nabla_{\boldsymbol{x}}f_{0}^{J}\left(\boldsymbol{x}^{*},\boldsymbol{\lambda}^{*}\right)+\sum_{i}\widetilde{\lambda}_{i}\nabla_{\boldsymbol{x}}f_{i}^{J}(\boldsymbol{x}^{*},\boldsymbol{\lambda}^{*}) & =\boldsymbol{0},\nonumber \\
\nabla_{\boldsymbol{\lambda}}f_{0}^{J}\left(\boldsymbol{x}^{*},\boldsymbol{\lambda}^{*}\right)+\sum_{i}\widetilde{\lambda}_{i}\nabla_{\boldsymbol{\lambda}}f_{i}^{J}(\boldsymbol{x}^{*},\boldsymbol{\lambda}^{*}) & =\boldsymbol{0},\nonumber \\
\widetilde{\lambda}_{i}f_{i}^{J}(\boldsymbol{x}^{*},\boldsymbol{\lambda}^{*}) & =0,\forall i\nonumber \\
f_{i}^{J}(\boldsymbol{x}^{*},\boldsymbol{\lambda}^{*}) & \leq0,\forall i,\label{eq:LKKT}
\end{align}
Then the primal-dual pair $\left(\boldsymbol{x}^{*}\in\mathcal{X},\Theta^{*}=\left\{ \boldsymbol{y}^{J}\left(\boldsymbol{x}^{*},\boldsymbol{\lambda}^{*},\boldsymbol{\xi}\right),\forall\boldsymbol{\xi}\right\} \right)$
and $\boldsymbol{\lambda}^{*}$ satisfies the KKT conditions in (\ref{eq:STscon}),
(\ref{eq:LTscon}) and (\ref{eq:CompS}) up to an error of $O\left(e\left(J\right)\right)$,
where $\lim_{J\rightarrow\infty}e\left(J\right)=0$, providing that
the following linear independence regularity condition (LIRC) holds:
The gradients of the long-term constraints $\nabla_{\boldsymbol{\lambda}}f_{i}^{J}(\boldsymbol{x}^{*},\boldsymbol{\lambda}^{*}),\forall i$
are linearly independent, and the gradients of the short-term constraints
$\partial_{\boldsymbol{y}}h_{j}\left(\boldsymbol{y}^{J}\left(\boldsymbol{x}^{*},\boldsymbol{\lambda}^{*},\boldsymbol{\xi}\right),\boldsymbol{\xi}\right),\forall j$
are also linearly independent for every $\boldsymbol{\xi}\in\Omega$
outside a set of probability zero. 
\end{thm}

Please refer to Appendix \ref{subsec:Proof-of-Lemma-surr} for the
proof.

The LIRC in Theorem \ref{thm:PDprac} is used to guarantee the complementary
slackness in (\ref{eq:CompS}). It is similar to another well known
regularity condition, the linear independence constraint qualification
(LICQ)\footnote{LICQ ensures the\textbf{ }existence of KKT point (i.e., stationary
point that satisfies KKT conditions) for a smooth optimization problem.}, in the sense that they both require linear independence conditions
on some gradients. However, the details are different, e.g., LICQ
requires that the gradients of the active inequality constraints and
the gradients of the equality constraints are linearly independent.
On the other hand, the condition that $\boldsymbol{y}^{J}\left(\boldsymbol{x},\boldsymbol{\lambda},\boldsymbol{\xi}\right)$
is continuously differentiable can be guaranteed by the short-term
sub-algorithm design, as will be detailed in Section \ref{subsec:short-term algorithm}. 

\section{Primal-Dual Decomposition based Stochastic Successive Convex Approximation\label{sec:Constrained-Stochastic-Successiv}}

Theorem \ref{thm:PDprac} states that a KKT solution (up to error
$e\left(J\right)$) can be found by first solving a stationary point
$\boldsymbol{x}^{*},\boldsymbol{\lambda}^{*}$ of the relaxed long-term
problem $\mathcal{P}_{L}^{J}$, and then finding a stationary point
$\boldsymbol{y}^{J}\left(\boldsymbol{x},\boldsymbol{\lambda},\boldsymbol{\xi}\right)$
of $\mathcal{P}_{S}\left(\boldsymbol{x}^{*},\boldsymbol{\lambda}^{*},\boldsymbol{\xi}\right)$
for each $\boldsymbol{\xi}$. In this section, we propose PDD-SSCA
to find KKT solutions of $\mathcal{P}$. Specifically, PDD-SSCA contains
two sub-algorithms: a long-term sub-algorithm and a short-term sub-algorithm.
We first present the long-term sub-algorithm which converges to a
stationary point of the non-convex constrained stochastic optimization
problem $\mathcal{P}_{L}^{J}$. Then, we discuss several general methods
to design the short-term sub-algorithm that finds a stationary point
of $\mathcal{P}_{S}\left(\boldsymbol{x},\boldsymbol{\lambda},\boldsymbol{\xi}\right)$.
Finally, we establish the convergence of the overall algorithm. The
implementation details for the overall algorithm are provided in Section
\ref{sec:Deep-Neural-Network}.

\subsection{Long-term Sub-Algorithm for $\mathcal{P}_{L}^{J}$}

The long-term problem $\mathcal{P}_{L}^{J}$ is a single-stage stochastic
optimization problem with non-convex constraints. The conventional
single-stage stochastic SCA algorithms in \cite{NIPS2013_5129_StochasticMM,Yang_TSP2016_SSCA}
only consider deterministic and convex constraints. Recently, a constrained
stochastic successive convex approximation (CSSCA) framework is proposed
in \cite{Liu_TSP2017_CSSCA} to find a stationary point for a non-convex
constrained single-stage stochastic optimization problem. The long-term
sub-algorithm in this paper is based on the CSSCA framework and is
summarized in Algorithm 1. In the $t$-th iteration, the long-term
variables $\boldsymbol{x},\boldsymbol{\lambda}$ are updated by solving
a convex optimization problem obtained by replacing the objective
and constraint functions $f_{i}^{J}\left(\boldsymbol{x},\boldsymbol{\lambda}\right),i=0,1,...,m$
with their convex surrogate functions $\bar{f}_{i}^{t}\left(\boldsymbol{x},\boldsymbol{\lambda}\right),i=0,1,...,m$,
as elaborated below. 

In Step 1 of the $t$-th iteration, one random mini-batch $\left\{ \boldsymbol{\xi}_{j}^{t},j=1,...,B\right\} $
of $B$ state samples are obtained and the surrogate functions $\bar{f}_{i}^{t}\left(\boldsymbol{x},\boldsymbol{\lambda}\right)$
are constructed based on the mini-batch $\left\{ \boldsymbol{\xi}_{j}^{t}\right\} $
and the current iterate $\boldsymbol{x}^{t},\boldsymbol{\lambda}^{t}$.
The surrogate functions $\bar{f}_{i}^{t}\left(\boldsymbol{x},\boldsymbol{\lambda}\right)$
can be viewed as convex approximations of the objective and constraint
functions $f_{i}^{J}(\boldsymbol{x},\boldsymbol{\lambda}),\forall i$
of the long-term problem $\mathcal{P}_{L}^{J}$. Specifically, divide
$g_{i}\left(\boldsymbol{x},\boldsymbol{y},\boldsymbol{\xi}\right)$
into two components as
\[
g_{i}\left(\boldsymbol{x},\boldsymbol{y},\boldsymbol{\xi}\right)=g_{i}^{c}\left(\boldsymbol{x},\boldsymbol{y},\xi\right)+g_{i}^{\bar{c}}\left(\boldsymbol{x},\boldsymbol{y},\xi\right),
\]
such that $g_{i}^{c}\left(\boldsymbol{x},\boldsymbol{y},\xi\right)$
is convex w.r.t. $\boldsymbol{x}$ and $g_{i}^{\bar{c}}\left(\boldsymbol{x},\boldsymbol{y},\xi\right)$
can be either convex or non-convex. Then a structured surrogate function
$\bar{f}_{i}^{t}\left(\boldsymbol{x},\boldsymbol{\lambda}\right)$
can be constructed as \cite{Liu_TSP2017_CSSCA}
\begin{align}
\bar{f}_{i}^{t}\left(\boldsymbol{x},\boldsymbol{\lambda}\right) & =\left(1-\rho^{t}\right)f_{i}^{t-1}+\rho^{t}\frac{1}{B}\sum_{j=1}^{B}\bigg[g_{i}^{c}\left(\boldsymbol{x},\boldsymbol{y}_{j}^{t},\boldsymbol{\xi}_{j}^{t}\right)\nonumber \\
 & +g_{i}^{\bar{c}}\left(\boldsymbol{x}^{t},\boldsymbol{y}_{j}^{t},\boldsymbol{\xi}_{j}^{t}\right)+\partial_{\boldsymbol{x}}^{T}g_{i}^{\bar{c}}\left(\boldsymbol{x}^{t},\boldsymbol{y}_{j}^{t},\boldsymbol{\xi}_{j}^{t}\right)\left(\boldsymbol{x}-\boldsymbol{x}^{t}\right)\bigg]\nonumber \\
 & +\left(\left(1-\rho^{t}\right)\mathbf{f}_{x,i}^{t-1}+\mathbf{f}_{y,i}^{t}\right)^{T}\left(\boldsymbol{x}-\boldsymbol{x}^{t}\right)\nonumber \\
 & +(\mathbf{f}_{\lambda,i}^{t})^{T}(\boldsymbol{\lambda}-\boldsymbol{\lambda}^{t})+\tau_{i}\left(\left\Vert \boldsymbol{x}-\boldsymbol{x}^{t}\right\Vert ^{2}+\left\Vert \boldsymbol{\lambda}-\boldsymbol{\lambda}^{t}\right\Vert ^{2}\right),\label{eq:SSF}
\end{align}
where $\boldsymbol{y}_{j}^{t}$ is an abbreviation for $\boldsymbol{y}^{J}\left(\boldsymbol{x}^{t},\boldsymbol{\lambda}^{t},\boldsymbol{\xi}_{j}^{t}\right)$,
$\left\{ \rho^{t}\in\left(0,1\right]\right\} $ is a decreasing sequence
satisfying $\rho^{t}\rightarrow0$, $\sum_{t}\rho^{t}=\infty$, $\sum_{t}\left(\rho^{t}\right)^{2}<\infty$,
$\tau_{i}>0$ can be any constant, $f_{i}^{t}$ is an approximation
for $f_{i}^{J}(\boldsymbol{x}^{t},\boldsymbol{\lambda}^{t})$ and
it is updated recursively according to
\[
f_{i}^{t}=\left(1-\rho^{t}\right)f_{i}^{t-1}+\rho^{t}\frac{1}{B}\sum_{j=1}^{B}g_{i}\left(\boldsymbol{x}^{t},\boldsymbol{y}_{j}^{t},\boldsymbol{\xi}_{j}^{t}\right),
\]
with $f_{i}^{-1}=0$, and $\mathbf{f}_{x,i}^{t},\mathbf{f}_{y,i}^{t},\mathbf{f}_{\lambda,i}^{t}$
are approximations for the gradients $\mathbb{E}\left[\partial_{\boldsymbol{x}}g_{i}\left(\boldsymbol{x},\boldsymbol{y}^{J}\left(\boldsymbol{x},\boldsymbol{\lambda},\boldsymbol{\xi}\right),\boldsymbol{\xi}\right)\right]$,
$\mathbb{E}\left[\partial_{\boldsymbol{x}}\boldsymbol{y}^{J}\left(\boldsymbol{x},\boldsymbol{\lambda},\boldsymbol{\xi}\right)\partial_{\boldsymbol{y}}g_{i}\left(\boldsymbol{x},\boldsymbol{y}^{J}\left(\boldsymbol{x},\boldsymbol{\lambda},\boldsymbol{\xi}\right),\boldsymbol{\xi}\right)\right]$
and $\mathbb{E}\left[\nabla_{\boldsymbol{\lambda}}g_{i}\left(\boldsymbol{x},\boldsymbol{y}^{J}\left(\boldsymbol{x},\boldsymbol{\lambda},\boldsymbol{\xi}\right),\boldsymbol{\xi}\right)\right]$,
respectively, which are updated recursively according to
\begin{align}
\mathbf{f}_{x,i}^{t} & =\left(1-\rho^{t}\right)\mathbf{f}_{x,i}^{t-1}+\rho^{t}\frac{1}{B}\sum_{j=1}^{B}\partial_{\boldsymbol{x}}g_{i}\left(\boldsymbol{x}^{t},\boldsymbol{y}_{j}^{t},\boldsymbol{\xi}_{j}^{t}\right),\nonumber \\
\mathbf{f}_{y,i}^{t} & =\left(1-\rho^{t}\right)\mathbf{f}_{y,i}^{t-1}+\rho^{t}\frac{1}{B}\sum_{j=1}^{B}\partial_{\boldsymbol{x}}\boldsymbol{y}_{j}^{t}\partial_{\boldsymbol{y}}g_{i}\left(\boldsymbol{x}^{t},\boldsymbol{y}_{j}^{t},\boldsymbol{\xi}_{j}^{t}\right),\nonumber \\
\mathbf{f}_{\lambda,i}^{t} & =\left(1-\rho^{t}\right)\mathbf{f}_{\lambda,i}^{t-1}+\rho^{t}\frac{1}{B}\sum_{j=1}^{B}\partial_{\boldsymbol{\lambda}}\boldsymbol{y}_{j}^{t}\partial_{\boldsymbol{y}}g_{i}\left(\boldsymbol{x}^{t},\boldsymbol{y}_{j}^{t},\boldsymbol{\xi}_{j}^{t}\right),\label{eq:fyx}
\end{align}
with $\mathbf{f}_{i}^{-1}=\boldsymbol{0}$, $\partial_{\boldsymbol{x}}\boldsymbol{y}_{j}^{t}=\partial_{\boldsymbol{x}}\boldsymbol{y}^{J}\left(\boldsymbol{x}^{t},\boldsymbol{\lambda}^{t},\boldsymbol{\xi}_{j}^{t}\right)\in\mathbb{R}^{n_{x}\times n_{y}}$
($\partial_{\boldsymbol{\lambda}}\boldsymbol{y}_{j}^{t}=\partial_{\boldsymbol{\lambda}}\boldsymbol{y}^{J}\left(\boldsymbol{x}^{t},\boldsymbol{\lambda}^{t},\boldsymbol{\xi}_{j}^{t}\right)\in\mathbb{R}^{m\times n_{y}}$)
is the derivative of the vector function $\boldsymbol{y}^{J}\left(\boldsymbol{x}^{t},\boldsymbol{\lambda}^{t},\boldsymbol{\xi}_{j}^{t}\right)$
to the vector $\boldsymbol{x}$ ($\boldsymbol{\lambda}$). Later in
Lemma \ref{lem:Property-surrogate}, we will show that $\lim_{t\rightarrow\infty}f_{i}^{t}-f_{i}^{J}(\boldsymbol{x}^{t},\boldsymbol{\lambda}^{t})=0$,
$\lim_{t\rightarrow\infty}\mathbf{f}_{x,i}^{t}+\mathbf{f}_{y,i}^{t}-\nabla_{\boldsymbol{x}}f_{i}^{J}(\boldsymbol{x}^{t},\boldsymbol{\lambda}^{t})=\boldsymbol{0}$
and $\lim_{t\rightarrow\infty}\mathbf{f}_{\lambda,i}^{t}-\nabla_{\boldsymbol{\lambda}}f_{i}^{J}(\boldsymbol{x}^{t},\boldsymbol{\lambda}^{t})=\boldsymbol{0}$.
Therefore, the function value and gradient of $\bar{f}_{i}^{t}\left(\boldsymbol{x},\boldsymbol{\lambda}\right)$
is consistent with that of the original function $f_{i}^{J}(\boldsymbol{x},\boldsymbol{\lambda})$
at the current iterate $\boldsymbol{x}^{t},\boldsymbol{\lambda}^{t}$,
which is the key to guarantee the convergence of the algorithm to
a stationary point. The structured surrogate function in (\ref{eq:SSF})
contains the convex component $g_{i}^{c}\left(\boldsymbol{x},\boldsymbol{y}_{j}^{t},\boldsymbol{\xi}_{j}^{t}\right)$
of the original sample objective function $g_{i}\left(\boldsymbol{x},\boldsymbol{y}_{j}^{t},\boldsymbol{\xi}_{j}^{t}\right)$,
which may help to reduce the approximation error and potentially achieve
a faster initial convergence speed \cite{Yang_TSP2016_SSCA}. How
to divide $g_{i}\left(\boldsymbol{x},\boldsymbol{y},\boldsymbol{\xi}\right)$
into two components to achieve a good initial convergence speed depends
on the specific problem.

In Step 2 of the $t$-th iteration, the optimal solution $\bar{\boldsymbol{x}}^{t},\bar{\boldsymbol{\lambda}}^{t}$
of the following problem is solved:
\begin{align}
\left(\bar{\boldsymbol{x}}^{t},\bar{\boldsymbol{\lambda}}^{t}\right)=\underset{\boldsymbol{x}\in\mathcal{X},\boldsymbol{\lambda}\succeq\boldsymbol{0}}{\text{argmin}}\: & \bar{f}_{0}^{t}\left(\boldsymbol{x},\boldsymbol{\lambda}\right)\label{eq:Pitert}\\
s.t.\: & \bar{f}_{i}^{t}\left(\boldsymbol{x},\boldsymbol{\lambda}\right)\leq0,\forall i,\nonumber 
\end{align}
which is a convex approximation of $\mathcal{P}_{L}^{J}$. Note that
Problem (\ref{eq:Pitert}) is not necessarily feasible. If Problem
(\ref{eq:Pitert}) turns out to be infeasible, the optimal solution
$\left(\bar{\boldsymbol{x}}^{t},\bar{\boldsymbol{\lambda}}^{t}\right)$
of the following convex problem is solved: 

\begin{align}
\left(\bar{\boldsymbol{x}}^{t},\bar{\boldsymbol{\lambda}}^{t}\right)=\underset{\boldsymbol{x}\in\mathcal{X},\boldsymbol{\lambda}\succeq\boldsymbol{0},\alpha}{\text{argmin}} & \:\alpha\label{eq:Pitert-1}\\
s.t.\: & \bar{f}_{i}^{t}\left(\boldsymbol{x}\right)\leq\alpha,\forall i,\nonumber 
\end{align}
which minimizes the approximate constraint functions. 

Finally, in Step 3, given $\bar{\boldsymbol{x}}^{t},\bar{\boldsymbol{\lambda}}^{t}$
in one of the above two cases, $\boldsymbol{x},\boldsymbol{\lambda}$
is updated according to
\begin{align}
\boldsymbol{x}^{t+1} & =\left(1-\gamma^{t}\right)\boldsymbol{x}^{t}+\gamma^{t}\bar{\boldsymbol{x}}^{t}.\nonumber \\
\boldsymbol{\lambda}^{t+1} & =\left(1-\gamma^{t}\right)\boldsymbol{\lambda}^{t}+\gamma^{t}\bar{\boldsymbol{\lambda}}^{t}.\label{eq:updatext}
\end{align}
where $\left\{ \gamma^{t}\in\left(0,1\right]\right\} $ is a decreasing
sequence satisfying $\gamma^{t}\rightarrow0$, $\sum_{t}\gamma^{t}=\infty$,
$\sum_{t}\left(\gamma^{t}\right)^{2}<\infty$, $\lim_{t\rightarrow\infty}\gamma^{t}/\rho^{t}=0$. 

\begin{algorithm}
\caption{\label{alg1}Long-term Sub-Algorithm for $\mathcal{P}_{L}^{J}$}

\textbf{Input: }$\left\{ \rho^{t},\gamma^{t}\right\} $.

\textbf{Initialize:} $\boldsymbol{x}^{0}\in\mathcal{X}$,$\boldsymbol{\lambda}^{0}\succeq\boldsymbol{0}$;
$t=0$.

\textbf{Step 1: }

Obtain a mini-batch $\left\{ \boldsymbol{\xi}_{j}^{t},j=1,...,B\right\} $. 

Construct\textbf{ }the surrogate functions $\bar{f}_{i}^{t}\left(\boldsymbol{x},\boldsymbol{\lambda}\right),\forall i$
according to (\ref{eq:SSF}).

\textbf{Step 2: }

//Objective update

\textbf{If} Problem (\ref{eq:Pitert}) is feasible

Solve (\ref{eq:Pitert}) to obtain $\bar{\boldsymbol{x}}^{t},\bar{\boldsymbol{\lambda}}^{t}$.

//Feasible update

\textbf{Else }

Solve (\ref{eq:Pitert-1}) to obtain $\bar{\boldsymbol{x}}^{t},\bar{\boldsymbol{\lambda}}^{t}$.

\textbf{End if}

\textbf{Step 3: }

Update $\boldsymbol{x}^{t+1},\boldsymbol{\lambda}^{t+1}$ according
to (\ref{eq:updatext}).

\textbf{Step 4: }

\textbf{Let} $t=t+1$ and return to Step 1.
\end{algorithm}

\subsection{Short-term Sub-Algorithm for $\mathcal{P}_{S}$\label{subsec:short-term algorithm}}

The short-term subproblem $\mathcal{P}_{S}\left(\boldsymbol{x},\boldsymbol{\lambda},\boldsymbol{\xi}\right)$
can be solved using existing deterministic optimization algorithms.
In general, a deterministic and iterative short-term sub-algorithm
starts from an initial point $\boldsymbol{y}^{0}\left(\boldsymbol{x},\boldsymbol{\lambda},\boldsymbol{\xi}\right)$,
and then generates a sequence $\left\{ \boldsymbol{y}^{j}\left(\boldsymbol{x},\boldsymbol{\lambda},\boldsymbol{\xi}\right)\right\} $
of iterates that converge to a stationary point of the short-term
subproblem $\mathcal{P}_{S}\left(\boldsymbol{x},\boldsymbol{\lambda},\boldsymbol{\xi}\right)$.
Specifically, the $j$-th iteration of a general short-term sub-algorithm
can be expressed as a mapping from $\boldsymbol{y}^{j-1}\left(\boldsymbol{x},\boldsymbol{\lambda},\boldsymbol{\xi}\right)$
to $\boldsymbol{y}^{j}\left(\boldsymbol{x},\boldsymbol{\lambda},\boldsymbol{\xi}\right)$
as
\begin{equation}
\boldsymbol{y}^{j}\left(\boldsymbol{x},\boldsymbol{\lambda},\boldsymbol{\xi}\right)=\mathcal{A}^{j}\left(\boldsymbol{y}^{j-1}\left(\boldsymbol{x},\boldsymbol{\lambda},\boldsymbol{\xi}\right),\boldsymbol{x},\boldsymbol{\lambda},\boldsymbol{\xi}\right),j=1,2,...,\label{eq:Smap}
\end{equation}
which depends on problem parameters $\boldsymbol{x},\boldsymbol{\lambda},\boldsymbol{\xi}$. 

To avoid confusion with the iteration of the long-term sub-algorithm,
an iteration of the short-term sub-algorithm will be called an \textit{inner
iteration}. To ensure the convergence of the overall algorithm, the
short-term sub-algorithm is assumed to satisfy the following conditions.

\begin{assumption}[Assumptions on the short-term sub-algorithm]\label{asm:shortAlg}$\:$
\begin{enumerate}
\item The initial point $\boldsymbol{y}^{0}\left(\boldsymbol{x},\boldsymbol{\lambda},\boldsymbol{\xi}\right)$
is differentiable w.r.t. $\boldsymbol{x},\boldsymbol{\lambda}$, w.p.1.
\item For any $\boldsymbol{x}\in\mathcal{X},\boldsymbol{y}\in\mathcal{Y},\boldsymbol{\lambda}\succeq\boldsymbol{0}$
and iteration number $j$, $\mathcal{A}^{j}\left(\boldsymbol{y},\boldsymbol{x},\boldsymbol{\lambda},\boldsymbol{\xi}\right)$
is differentiable w.r.t. $\boldsymbol{y},\boldsymbol{x},\boldsymbol{\lambda}$,
w.p.1.
\item For any $\boldsymbol{x}\in\mathcal{X},\boldsymbol{\lambda}\succeq\boldsymbol{0}$,
the sequence $\left\{ \boldsymbol{y}^{j}\left(\boldsymbol{x},\boldsymbol{\lambda},\boldsymbol{\xi}\right)\right\} $
converges to a stationary set of $\mathcal{P}_{S}\left(\boldsymbol{x},\boldsymbol{\lambda},\boldsymbol{\xi}\right)$,
w.p.1.
\end{enumerate}
\end{assumption}

Assumption \ref{asm:shortAlg} ensures that $\boldsymbol{y}^{J}\left(\boldsymbol{x},\boldsymbol{\lambda},\boldsymbol{\xi}\right)$
is differentiable w.r.t. $\boldsymbol{x},\boldsymbol{\lambda}$ for
any finite $J$ and $\boldsymbol{y}^{J}\left(\boldsymbol{x},\boldsymbol{\lambda},\boldsymbol{\xi}\right)$
satisfies the KKT conditions of $\mathcal{P}_{S}\left(\boldsymbol{x},\boldsymbol{\lambda},\boldsymbol{\xi}\right)$
(up to certain tolerable error $O\left(e(J)\right)$) as in (\ref{eq:STsconeJ}).
The first condition in Assumption \ref{asm:shortAlg} can be satisfied
by a proper choice of the initial points $\boldsymbol{y}^{0}\left(\boldsymbol{x},\boldsymbol{\lambda},\boldsymbol{\xi}\right),\forall\boldsymbol{x},\boldsymbol{\lambda},\boldsymbol{\xi}$.
The second and third conditions are satisfied by many standard iterative
algorithms that are designed to find stationary points of a non-convex
problem. In the following, we give two examples of short-term sub-algorithms
that satisfy Assumption \ref{asm:shortAlg}.

\subsubsection{Gradient Projection Algorithm}

When the feasible set of the short-term subproblem $\mathcal{P}_{S}\left(\boldsymbol{x},\boldsymbol{\lambda},\boldsymbol{\xi}\right)$,
denoted by $\mathcal{Y}\left(\boldsymbol{x},\boldsymbol{\lambda},\boldsymbol{\xi}\right)\triangleq\left\{ \boldsymbol{y}\in\mathcal{Y}:h_{j}\left(\boldsymbol{y}\left(\boldsymbol{\xi}\right),\boldsymbol{\xi}\right)\leq0,\:\forall j\right\} $,
is convex, the gradient projection (GP) algorithm \cite{Bertsekas_book99_NProgramming}
can be used to find a stationary point of $\mathcal{P}_{S}\left(\boldsymbol{x},\boldsymbol{\lambda},\boldsymbol{\xi}\right)$.
Note that the objective function of $\mathcal{P}_{S}\left(\boldsymbol{x},\boldsymbol{\lambda},\boldsymbol{\xi}\right)$
can still be non-convex. When the GP algorithm is used as the short-term
sub-algorithm, the mapping $\mathcal{A}^{j}\left(\boldsymbol{y}^{j-1}\left(\boldsymbol{x},\boldsymbol{\lambda},\boldsymbol{\xi}\right),\boldsymbol{x},\boldsymbol{\lambda},\boldsymbol{\xi}\right)$
for the $j$-th inner iteration is given by
\begin{equation}
\mathcal{A}^{j}\left(\boldsymbol{y}^{j-1}\right)=\mathbb{P}_{\mathcal{Y}\left(\boldsymbol{x},\boldsymbol{\lambda},\boldsymbol{\xi}\right)}\left[\boldsymbol{y}^{j-1}-\alpha_{j}\partial_{\boldsymbol{y}}g_{s}\left(\boldsymbol{x},\boldsymbol{\lambda},\boldsymbol{y}^{j-1},\boldsymbol{\xi}\right)\right],\label{eq:GP}
\end{equation}
where $g_{s}\left(\boldsymbol{x},\boldsymbol{\lambda},\boldsymbol{y},\boldsymbol{\xi}\right)=g_{0}\left(\boldsymbol{x},\boldsymbol{y},\boldsymbol{\xi}\right)+\sum_{i}\lambda_{i}g_{i}\left(\boldsymbol{x},\boldsymbol{y},\boldsymbol{\xi}\right)$,
$\mathbb{P}_{\mathcal{Y}\left(\boldsymbol{x},\boldsymbol{\lambda},\boldsymbol{\xi}\right)}$
is the projection onto the feasible set $\mathcal{Y}\left(\boldsymbol{x},\boldsymbol{\lambda},\boldsymbol{\xi}\right)$,
and $\left\{ \alpha_{j}>0,j=1,...,J\right\} $ is a properly chosen
step size sequence, e.g., a diminishing step size sequence \cite{Bertsekas_book99_NProgramming}.
Note that we have omitted $\left(\boldsymbol{x},\boldsymbol{\lambda},\boldsymbol{\xi}\right)$
in the mapping $\mathcal{A}^{j}$ for simplicity of notation. It can
be verified that the GP algorithm satisfies Assumption \ref{asm:shortAlg}. 

\subsubsection{Majorization-Minimization Algorithm }

Majorization-Minimization (MM) \cite{Meisam_thesis14_BSUM,Sun_TSP2017_MM}
can be used to find a stationary point of $\mathcal{P}_{S}\left(\boldsymbol{x},\boldsymbol{\lambda},\boldsymbol{\xi}\right)$
for general cases. When MM is used as the short-term sub-algorithm,
the mapping $\mathcal{A}^{j}\left(\boldsymbol{y}^{j-1}\left(\boldsymbol{x},\boldsymbol{\xi}\right),\boldsymbol{x},\boldsymbol{\lambda},\boldsymbol{\xi}\right)$
is given by
\begin{align}
\mathcal{A}^{j}\left(\boldsymbol{y}^{j-1}\right)= & \underset{\boldsymbol{y}}{\text{argmin}\:}u_{s}\left(\boldsymbol{y};\boldsymbol{y}^{j-1},\boldsymbol{x},\boldsymbol{\lambda}\right)\label{eq:cvxMM}\\
 & \text{s.t. }u_{i}\left(\boldsymbol{y};\boldsymbol{y}^{j-1}\right)\leq0,\:i=1,...,n,\nonumber 
\end{align}
where $u_{s}\left(\boldsymbol{y};\boldsymbol{y}^{j-1},\boldsymbol{x},\boldsymbol{\lambda}\right)$
is a surrogate function of $g_{s}\left(\boldsymbol{x},\boldsymbol{\lambda},\boldsymbol{y},\boldsymbol{\xi}\right)$,
$u_{i}\left(\boldsymbol{y};\boldsymbol{y}^{j-1}\right)$ is a surrogate
function of $h_{i}\left(\boldsymbol{y},\boldsymbol{\xi}\right)$ for
$i=1,...,n$, satisfying the following conditions:
\begin{enumerate}
\item $u_{s}\left(\boldsymbol{y}^{j-1};\boldsymbol{y}^{j-1},\boldsymbol{x},\boldsymbol{\lambda}\right)=g_{s}\left(\boldsymbol{x},\boldsymbol{\lambda},\boldsymbol{y}^{j-1},\boldsymbol{\xi}\right)$,
$u_{s}\left(\boldsymbol{y};\boldsymbol{y}^{j-1},\boldsymbol{x},\boldsymbol{\lambda}\right)\geq g_{s}\left(\boldsymbol{x},\boldsymbol{\lambda},\boldsymbol{y},\boldsymbol{\xi}\right),\forall\boldsymbol{y}\in\mathcal{Y}$
and $\nabla u_{s}\left(\boldsymbol{y}^{j-1};\boldsymbol{y}^{j-1},\boldsymbol{x},\boldsymbol{\lambda}\right)=\partial_{\boldsymbol{y}}g_{i}\left(\boldsymbol{x},\boldsymbol{\lambda},\boldsymbol{y}^{j-1},\boldsymbol{\xi}\right)$.
\item $u_{i}\left(\boldsymbol{y}^{j-1};\boldsymbol{y}^{j-1}\right)=h_{i}\left(\boldsymbol{y}^{j-1},\boldsymbol{\xi}\right)$,
$u_{i}\left(\boldsymbol{y};\boldsymbol{y}^{j-1}\right)\geq h_{i}\left(\boldsymbol{y},\boldsymbol{\xi}\right),\forall\boldsymbol{y}\in\mathcal{Y}$
and $\nabla u_{i}\left(\boldsymbol{y}^{j-1};\boldsymbol{y}^{j-1}\right)=\partial_{\boldsymbol{y}}h_{i}\left(\boldsymbol{y}^{j-1},\boldsymbol{\xi}\right)$.
\item $u_{i}\left(\boldsymbol{y};\boldsymbol{y}^{j-1}\right)$ is uniformly
strongly convex in $\boldsymbol{y}$.
\item $u_{s}\left(\boldsymbol{y};\boldsymbol{y}^{'},\boldsymbol{x},\boldsymbol{\lambda}\right)$
and $u_{i}\left(\boldsymbol{y};\boldsymbol{y}^{'}\right),\forall i$
are differentiable w.r.t. $\boldsymbol{y}^{'},\boldsymbol{y},\boldsymbol{x},\boldsymbol{\lambda}$.
\end{enumerate}

A simple example surrogate function that satisfies the above four
conditions is
\begin{align*}
u_{s}\left(\boldsymbol{y};\boldsymbol{y}^{'},\boldsymbol{x},\boldsymbol{\lambda}\right) & =g_{s}\left(\boldsymbol{x},\boldsymbol{\lambda},\boldsymbol{y}^{'},\boldsymbol{\xi}\right)+\tau_{s}\left\Vert \boldsymbol{y}-\boldsymbol{y}^{'}\right\Vert ^{2}\\
 & +\partial_{\boldsymbol{y}}^{T}g_{s}\left(\boldsymbol{x},\boldsymbol{\lambda},\boldsymbol{y}^{'},\boldsymbol{\xi}\right)\left(\boldsymbol{y}-\boldsymbol{y}^{'}\right),
\end{align*}
\begin{align*}
u_{i}\left(\boldsymbol{y};\boldsymbol{y}^{'}\right) & =h_{i}\left(\boldsymbol{y}^{'},\boldsymbol{\xi}\right)+\partial_{\boldsymbol{y}}^{T}h_{i}\left(\boldsymbol{y}^{'},\boldsymbol{\xi}\right)\left(\boldsymbol{y}-\boldsymbol{y}^{'}\right)\\
 & +\tau_{s,i}\left\Vert \boldsymbol{y}-\boldsymbol{y}^{'}\right\Vert ^{2},
\end{align*}
where $\tau_{s},\tau_{s,i}>0$ is chosen to be a sufficiently large
number to satisfy the upper bound condition: $u_{s}\left(\boldsymbol{y};\boldsymbol{y}^{j-1},\boldsymbol{x},\boldsymbol{\lambda}\right)\geq g_{s}\left(\boldsymbol{x},\boldsymbol{\lambda},\boldsymbol{y},\boldsymbol{\xi}\right),\forall\boldsymbol{y}\in\mathcal{Y}$
and $u_{i}\left(\boldsymbol{y};\boldsymbol{y}^{j-1}\right)\geq h_{i}\left(\boldsymbol{y},\boldsymbol{\xi}\right),\forall\boldsymbol{y}\in\mathcal{Y}$.
Note that such $\tau_{s},\tau_{s,i}$ can always be found since the
second-order derivative of the functions $g_{s}\left(\boldsymbol{x},\boldsymbol{\lambda},\boldsymbol{y},\boldsymbol{\xi}\right)$
and $h_{i}\left(\boldsymbol{y},\boldsymbol{\xi}\right)$'s are assumed
to be uniformly bounded.

From the convergence result for the MM algorithm in \cite{Meisam_thesis14_BSUM}
and the above conditions on $u_{s}\left(\boldsymbol{y};\boldsymbol{y}^{j-1},\boldsymbol{x},\boldsymbol{\lambda}\right)$
and $u_{i}\left(\boldsymbol{y};\boldsymbol{y}^{j-1}\right),\forall i$,
one can verify that the MM algorithm satisfies Assumption \ref{asm:shortAlg}.

\subsubsection{Short-term Sub-Algorithm for Section \ref{subsec:Examples-of-Problem}\label{subsec:Short-term-Algorithm-for-Exm}}

In this subsection, we give concrete examples of the short-term sub-algorithm
for our applications in Section \ref{subsec:Examples-of-Problem}.

In Example 1, the short-term subproblem $\mathcal{P}_{S}\left(\boldsymbol{\lambda},\Upsilon,\boldsymbol{a},\boldsymbol{b}\right)$
is given by 
\begin{align}
\min & -\left[\log\left(1+\sum_{i=1}^{N}a_{i}p_{i}\right)\right]+\sum_{i=1}^{N}\lambda_{i}p_{i}+\Upsilon\sum_{i=1}^{N}b_{i}p_{i},\label{SP1}\\
\text{s.t.} & \Upsilon\geq0,\lambda_{i}\geq0,p_{i}\geq0,i=1,\ldots,N,\nonumber 
\end{align}
where $\lambda_{i}$ and $\Upsilon$ are the long-term Lagrange multipliers
associated with the long-term transmit power budget at SU $i$ and
the interference threshold constraint for the PU, respectively. The
optimal solution of $\mathcal{P}_{S}\left(\boldsymbol{\lambda},\Upsilon,\boldsymbol{a},\boldsymbol{b}\right)$
has a closed-form expression given by
\begin{equation}
p_{i}^{*}\left(\boldsymbol{\lambda},\Upsilon,\boldsymbol{a},\boldsymbol{b}\right)=\frac{1}{Na_{i}}\left(\frac{a_{i}}{b_{i}\Upsilon+\lambda_{i}}-1\right)^{+}.\label{eq:Sexam1}
\end{equation}

In Example 2, the short-term subproblem $\mathcal{P}_{S}\left(\boldsymbol{\theta},\boldsymbol{\lambda},\boldsymbol{H}\right)$
is given by

\begin{equation}
\underset{\boldsymbol{G}}{\text{min}}Tr\left(\boldsymbol{F}\boldsymbol{G}\boldsymbol{G}^{H}\boldsymbol{F}^{H}\right)-\sum_{k=1}^{K}\lambda_{k}r_{k}\left(\boldsymbol{\theta},\boldsymbol{G},\boldsymbol{H}\right),\label{eq:24}
\end{equation}
where $\boldsymbol{\lambda}=\left[\lambda_{1},...,\lambda_{K}\right]$
are the Lagrange multipliers associated with the average rate constraints.
A stationary point of $\mathcal{P}_{S}\left(\boldsymbol{\theta},\boldsymbol{\lambda},\boldsymbol{H}\right)$
can be found using the WMMSE algorithm \cite{Luo_TSP11_WMMSE}. The
basic idea is to first transform $\mathcal{P}_{S}\left(\boldsymbol{\theta},\boldsymbol{\lambda},\boldsymbol{H}\right)$
into the following WMMSE problem by introducing two auxiliary variables
$\boldsymbol{w},\boldsymbol{u}$:
\begin{equation}
\underset{\left\{ \boldsymbol{w},\boldsymbol{u},\boldsymbol{G}\right\} }{\textrm{min}}Tr\left(\boldsymbol{F}\boldsymbol{G}\boldsymbol{G}^{H}\boldsymbol{F}^{H}\right)+\sum_{k=1}^{K}\lambda_{k}(w_{k}e_{k}-\textrm{log}w_{k}),\label{eq:WMMSE}
\end{equation}
where $\boldsymbol{w}=\left[w_{1},...,w_{K}\right]^{T}$ with $w_{k}>0:\forall k$
is a weight vector for MSE; $\boldsymbol{u}=\left[u_{1},...,u_{K}\right]^{T}$
with $u_{k}$ denoting the receive coefficient; and

\begin{equation}
e_{k}\triangleq\left|u_{k}^{*}\boldsymbol{h}_{k}^{H}\boldsymbol{F}\boldsymbol{g}_{k}-1\right|^{2}+\sum_{i\neq k}\left|u_{k}^{*}\boldsymbol{h}_{k}^{H}\boldsymbol{F}\boldsymbol{g}_{i}\right|^{2}+\left|u_{k}\right|^{2},\label{eq:26}
\end{equation}
is the MSE of user $k.$ Problem (\ref{eq:WMMSE}) is convex in each
of the optimization variables $\boldsymbol{w},\boldsymbol{u},\boldsymbol{G}.$
Then, we can use the block coordinate descent method to solve (\ref{eq:WMMSE}).
Specifically, for given $\boldsymbol{G}$, the optimal $\boldsymbol{u}$
is given by the MMSE receive coefficient:

\begin{equation}
u_{k}=\left(\sum_{i=1}^{K}\left|\boldsymbol{h}_{k}^{H}\boldsymbol{F}\boldsymbol{g}_{i}\right|^{2}+1\right)^{-1}h_{k}^{H}\boldsymbol{F}\boldsymbol{g}_{k},\forall k.\label{eq:Sw}
\end{equation}
For given $\boldsymbol{G},\boldsymbol{u},$the optimal $\boldsymbol{w}$
is given by 
\begin{equation}
w_{k}=\left(1-u_{k}^{*}\boldsymbol{h}_{k}^{H}\boldsymbol{F}\boldsymbol{g}_{k}\right)^{-1},\forall k.\label{eq:sw}
\end{equation}
Finally, for given $\boldsymbol{w},\boldsymbol{u},$ the optimal $\boldsymbol{G}=\left[\boldsymbol{g}_{1},...,\boldsymbol{g}_{K}\right]$
is given by

\begin{equation}
\mathbf{g}_{k}=\left(\lambda_{k}w_{k}\sum_{i=1}^{K}\left|u_{i}\right|^{2}\boldsymbol{F}^{H}\boldsymbol{h}_{i}\boldsymbol{h}_{i}^{H}\boldsymbol{F}+\boldsymbol{I}\right)^{-1}\lambda_{k}w_{k}u_{k}\boldsymbol{F}^{H}\boldsymbol{h}_{k}.\label{eq:sG}
\end{equation}
After running (\ref{eq:Sw}) to (\ref{eq:sG}) for $J$ iterations,
we obtain an approximate stationary point of $\mathcal{P}_{S}\left(\boldsymbol{\theta},\boldsymbol{\lambda},\boldsymbol{H}\right)$
denoted as $\boldsymbol{G}^{J}\left(\boldsymbol{\theta},\boldsymbol{\lambda},\boldsymbol{H}\right)$.

\subsection{Overall Algorithm and Convergence Analysis}

In this section, we first describe the overall PDD-SSCA algorithm.
Then we analyze the convergence of the PDD-SSCA algorithm. The PDD-SSCA
algorithm first runs the long-term sub-algorithm to find a stationary
point of the long-term sub-problem $\boldsymbol{x}^{*},\boldsymbol{\lambda}^{*}$,
and then runs the short-term sub-algorithm to find a stationary point
$\boldsymbol{y}^{J}\left(\boldsymbol{x}^{*},\boldsymbol{\lambda}^{*},\boldsymbol{\xi}\right)$
of the short-term problem $\mathcal{P}_{S}\left(\boldsymbol{x}^{*},\boldsymbol{\lambda}^{*},\boldsymbol{\xi}\right)$
for each state realization $\boldsymbol{\xi}$ with fixed $\boldsymbol{x}^{*},\boldsymbol{\lambda}^{*}$.
In the following, we will show that the solution $\left(\boldsymbol{x}^{*},\Theta^{*}=\left\{ \boldsymbol{y}^{J}\left(\boldsymbol{x}^{*},\boldsymbol{\lambda}^{*},\boldsymbol{\xi}\right),\forall\boldsymbol{\xi}\right\} \right)$
found by the PDD-SSCA algorithm is a KKT solution of the original
two-stage stochastic optimization problem $\mathcal{P}$.

We first prove a key Lemma which establishes several important properties
of the surrogate functions $\bar{f}_{i}^{t}\left(\boldsymbol{x},\boldsymbol{\lambda}\right)$'s.
\begin{lem}
[Properties of surrogate functions]\label{lem:Property-surrogate}For
all $i\in\left\{ 0,...,m\right\} $ and $t=0,1,....$, we have
\begin{enumerate}
\item $\bar{f}_{i}^{t}\left(\boldsymbol{x},\boldsymbol{\lambda}\right)$
is uniformly strongly convex in $\boldsymbol{x},\boldsymbol{\lambda}$.
\item $\bar{f}_{i}^{t}\left(\boldsymbol{x},\boldsymbol{\lambda}\right)$
is a Lipschitz continuous function w.r.t. $\boldsymbol{x},\boldsymbol{\lambda}$.
Moreover, $\limsup_{t_{1},t_{2}\rightarrow\infty}\left|\bar{f}_{i}^{t_{1}}\left(\boldsymbol{x},\boldsymbol{\lambda}\right)-\bar{f}_{i}^{t_{2}}\left(\boldsymbol{x},\boldsymbol{\lambda}\right)\right|-B\sqrt{\left\Vert \boldsymbol{x}^{t_{1}}-\boldsymbol{x}^{t_{2}}\right\Vert ^{2}+\left\Vert \boldsymbol{\lambda}^{t_{1}}-\boldsymbol{\lambda}^{t_{2}}\right\Vert ^{2}}\leq0,\forall\boldsymbol{x}\in\mathcal{X},\boldsymbol{\lambda}\succeq\boldsymbol{0}$
for some constant $B>0$.
\item For any $\boldsymbol{x}\in\mathcal{X},\boldsymbol{\lambda}\succeq\boldsymbol{0}$,
the function $\bar{f}_{i}^{t}\left(\boldsymbol{x},\boldsymbol{\lambda}\right)$,
its derivative, and its second order derivative are uniformly bounded.
\item $\lim_{t\rightarrow\infty}\left|\bar{f}_{i}^{t}\left(\boldsymbol{x}^{t},\boldsymbol{\lambda}^{t}\right)-f_{i}^{J}\left(\boldsymbol{x}^{t},\boldsymbol{\lambda}^{t}\right)\right|=0$,
$\lim_{t\rightarrow\infty}\left\Vert \nabla_{\boldsymbol{x}}\bar{f}_{i}^{t}\left(\boldsymbol{x}^{t},\boldsymbol{\lambda}^{t}\right)-\nabla_{\boldsymbol{x}}f_{i}^{J}\left(\boldsymbol{x}^{t},\boldsymbol{\lambda}^{t}\right)\right\Vert =0$
and $\lim_{t\rightarrow\infty}\left\Vert \nabla_{\boldsymbol{\lambda}}\bar{f}_{i}^{t}\left(\boldsymbol{x}^{t},\boldsymbol{\lambda}^{t}\right)-\nabla_{\boldsymbol{\lambda}}f_{i}^{J}\left(\boldsymbol{x}^{t},\boldsymbol{\lambda}^{t}\right)\right\Vert =0$.
\item Consider a subsequence $\left\{ \boldsymbol{x}^{t_{j}},\boldsymbol{\lambda}^{t_{j}}\right\} _{j=1}^{\infty}$
converging to a limit point $\left(\boldsymbol{x}^{*},\boldsymbol{\lambda}^{*}\right)$.
There exist uniformly differentiable functions $\hat{f}_{i}\left(\boldsymbol{x},\boldsymbol{\lambda}\right)$
such that
\begin{align}
\lim_{j\rightarrow\infty}\bar{f}_{i}^{t_{j}}\left(\boldsymbol{x},\boldsymbol{\lambda}\right) & =\hat{f}_{i}\left(\boldsymbol{x},\boldsymbol{\lambda}\right),\:\forall\boldsymbol{x}\in\mathcal{X},\boldsymbol{\lambda}\succeq\boldsymbol{0},\label{eq:ghfhead}
\end{align}
almost surely.
\end{enumerate}
\end{lem}

Please refer to Appendix \ref{lem:Property-surrogate} for the proof.

The convergence analysis also relies on the Slater condition defined
below. 

\textbf{Slater condition:} Given a subsequence $\left\{ \boldsymbol{x}^{t_{j}},\boldsymbol{\lambda}^{t_{j}}\right\} _{j=1}^{\infty}$
converging to a limit point $\left(\boldsymbol{x}^{*},\boldsymbol{\lambda}^{*}\right)$
and let $\hat{f}_{i}\left(\boldsymbol{x},\boldsymbol{\lambda}\right),\forall i$
be the converged surrogate functions as defined in Lemma \ref{lem:Property-surrogate}.
We say that the Slater condition is satisfied at $\left(\boldsymbol{x}^{*},\boldsymbol{\lambda}^{*}\right)$
if there exists $\left(\boldsymbol{x}^{*},\boldsymbol{\lambda}^{*}\right)\in\textrm{relint}\mathcal{X}\times\mathbb{R}^{++}$
such that
\[
\hat{f}_{i}\left(\boldsymbol{x}^{*},\boldsymbol{\lambda}^{*}\right)<0,\:\forall i=1,...,m.
\]

With Lemma \ref{lem:Property-surrogate} and the Slater condition,
we are ready to prove the following main convergence result.
\begin{thm}
[Convergence of Algorithm 1]\label{thm:Convergence-of-Algorithm1-1}Suppose
Assumptions \ref{asm:convP} and \ref{asm:shortAlg} are satisfied,
and the initial point $\boldsymbol{x}^{0}\in\mathcal{X},\boldsymbol{\lambda}^{0}\in\mathbb{R}^{++}$
is a feasible point, i.e., $\max_{i\in\left\{ 1,...,m\right\} }f_{i}^{J}\left(\boldsymbol{x}^{0},\boldsymbol{\lambda}^{0}\right)\leq0$.
Let $\left\{ \boldsymbol{x}^{t},\boldsymbol{\lambda}^{t}\right\} _{t=1}^{\infty}$
denote the iterates generated by Algorithm 1 with a sufficiently small
initial step size $\gamma^{0}$. Then every limiting point $\left(\boldsymbol{x}^{*},\boldsymbol{\lambda}^{*}\right)$
of $\left\{ \boldsymbol{x}^{t},\boldsymbol{\lambda}^{t}\right\} _{t=1}^{\infty}$
satisfying the LIRC in Theorem \ref{thm:PDprac} and the Slater condition,
almost surely satisfies the KKT conditions in (\ref{eq:STscon}),
(\ref{eq:LTscon}) and (\ref{eq:CompS}) up to an error of $O\left(e\left(J\right)\right)$,
where $\lim_{J\rightarrow\infty}e\left(J\right)\rightarrow0$.
\end{thm}

\begin{IEEEproof}
It follows from Lemma \ref{lem:Property-surrogate} and the convergence
theorem of CSSCA in \cite{Liu_TSP2017_CSSCA} that, starting from
a feasible initial point, every limiting point $\left(\boldsymbol{x}^{*},\boldsymbol{\lambda}^{*}\right)$
generated by Algorithm 1 is a stationary point $\boldsymbol{x}^{*},\boldsymbol{\lambda}^{*}$
of the long-term sub-problem $\mathcal{P}_{L}^{J}$ almost surely,
providing that the initial step size\textit{ $\gamma^{0}$} is sufficiently
small, and the Slater condition is satisfied for $\boldsymbol{x}^{*},\boldsymbol{\lambda}^{*}$.
Then, it follows from the relaxed primal-dual decomposition method
in Theorem \ref{thm:PDprac} that the solution $\left(\boldsymbol{x}^{*},\Theta^{*}=\left\{ \boldsymbol{y}^{J}\left(\boldsymbol{x}^{*},\boldsymbol{\lambda}^{*},\boldsymbol{\xi}\right),\forall\boldsymbol{\xi}\right\} \right)$
found by the PDD-SSCA (Algorithm 1) is a KKT solution of the original
problem $\mathcal{P}$, up to certain error $e\left(J\right)$ that
diminishes to zero as $J\rightarrow\infty$. 
\end{IEEEproof}

Note that due to the stochastic nature of the problem/algorithm, we
need to assume that the step size is sufficiently small to make it
easier to handle the randomness caused by the random state for tractable
convergence analysis and rigorous convergence proof. However, choosing
a small $\gamma^{0}$ is usually not mandatory for the practical convergence
of Algorithm 1 as it is a sufficient not a necessary condition. In
the simulations, we find that the algorithm can still converge even
when the initial step size $\gamma^{0}$ is not small. In fact, in
practice, we may prefer to choose a not very small $\gamma^{0}$ to
achieve a faster initial convergence speed. Finally, we discuss the
convergence behavior of Algorithm 1 with an infeasible initial point.
Again, starting from a feasible initial point is just a sufficient
not a necessary condition. Due to the feasible update in (\ref{eq:Pitert-1}),
Algorithm 1 still converges to a KKT solution of Problem (\ref{eq:mainP})
with high probability, even when the initial point is infeasible \cite{Liu_TSP2017_CSSCA}.

\section{Deep Unrolling based Implementation\label{sec:Deep-Neural-Network}}

\subsection{Motivation of Deep Unrolling based Implementation}

In the long-term sub-algorithm, we need to calculate the gradient
of $\boldsymbol{y}^{J}\left(\boldsymbol{x},\boldsymbol{\lambda},\boldsymbol{\xi}\right)$
w.r.t. the long-term variables $\boldsymbol{x}$ and $\boldsymbol{\lambda}$:
$\partial_{\boldsymbol{x}}\boldsymbol{y}^{J}\left(\boldsymbol{x},\boldsymbol{\lambda},\boldsymbol{\xi}\right)$
and $\partial_{\boldsymbol{\lambda}}\boldsymbol{y}^{J}\left(\boldsymbol{x},\boldsymbol{\lambda},\boldsymbol{\xi}\right)$,
where $\boldsymbol{y}^{J}\left(\boldsymbol{x},\boldsymbol{\lambda},\boldsymbol{\xi}\right)$
is obtained by running the short-term sub-algorithm for $J$ iterations.
Since $\boldsymbol{y}^{J}\left(\boldsymbol{x},\boldsymbol{\lambda},\boldsymbol{\xi}\right)$
involves an iterative algorithm, it is usually not easy to calculate
its gradient in closed-form. One possible solution is to treat the
iterative short-term sub-algorithm as a black-box, and learn the mapping
$\boldsymbol{y}^{J}\left(\boldsymbol{x},\boldsymbol{\lambda},\boldsymbol{\xi}\right)$
between the input and the output by employing the DNN. Once we obtain
a DNN representation of the mapping $\boldsymbol{y}^{J}\left(\boldsymbol{x},\boldsymbol{\lambda},\boldsymbol{\xi}\right)$,
we can calculate the gradients $\partial_{\boldsymbol{x}}\boldsymbol{y}^{J}\left(\boldsymbol{x},\boldsymbol{\lambda},\boldsymbol{\xi}\right)$
and $\partial_{\boldsymbol{\lambda}}\boldsymbol{y}^{J}\left(\boldsymbol{x},\boldsymbol{\lambda},\boldsymbol{\xi}\right)$
using the well-known back propagation (BP) approach. Some representative
studies of such a solution can be found in \cite{Hong_TSP2018_WMMSEDNN,Yu_JSAC2019_Deepscheduling,Debbah_TSP2020_DLpowcont}
for different applications. For example, in \cite{Hong_TSP2018_WMMSEDNN},
the authors applied the multi-layer perceptron (MLP) and convolutional
neural network (CNN) to approximate the iterative WMMSE algorithm
used in Example 2. 

However, the black-box based DNNs suffer from poor interpretability
and generalization ability, and have no performance guarantee. Moreover,
the black-box based DNN often has a large number of parameters and
requires a lot of training samples, which incurs high training complexity
and memory overhead. To overcome such drawbacks, a number of works
\cite{Zhang_TSP19_Deepunrollpow,Li_TCI20_deepunroll,Hu_TWC20_DeepunrollWMMSE}
have proposed to unfold the iterations into a layer-wise structure
analogous to a NN based on the existing iterative algorithms. This
method is referred to as deep unrolling/unfolding and has a wide range
of applications in communications and signal processing. Compared
to the black-box based DNN, the deep unrolling based NN tends to have
better interpretability and generalization ability, as well as much
lower training complexity and memory overhead. Therefore, in this
paper, we propose to use the deep unrolling method to obtain a NN
representation of the mapping $\boldsymbol{y}^{J}\left(\boldsymbol{x},\boldsymbol{\lambda},\boldsymbol{\xi}\right)$
(or equivalently, the short-term sub-algorithm). 

\subsection{Deep Unrolling of the Short-term Sub-Algorithm \label{subsec:DNN-unroll} }

\begin{figure}
\begin{centering}
\includegraphics[width=85mm]{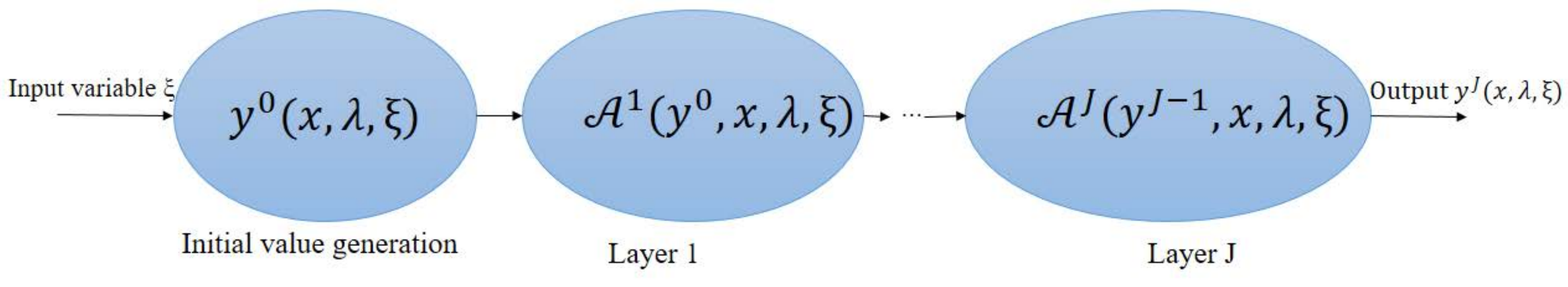}
\par\end{centering}
\caption{\label{fig:NNS}The NN representation of the short-term sub-algorithm.}
\end{figure}

We can apply the deep unrolling technique to obtain a NN representation
of the short-term sub-algorithm, where the input is $\boldsymbol{\xi}$,
the output is $\boldsymbol{y}^{J}\left(\boldsymbol{x},\boldsymbol{\lambda},\boldsymbol{\xi}\right)$,
the $0$-th layer is used to generate the initial value $\boldsymbol{y}^{0}\left(\boldsymbol{x},\boldsymbol{\lambda},\boldsymbol{\xi}\right)$,
and the $j$-th layer for $j=1,...,J$ is simply given by $\mathcal{A}^{j}\left(\boldsymbol{y}^{j-1}\left(\boldsymbol{x},\boldsymbol{\lambda},\boldsymbol{\xi}\right),\boldsymbol{x},\boldsymbol{\lambda},\boldsymbol{\xi}\right)$,
as illustrated in Fig. \ref{fig:NNS}.

\begin{figure}
\begin{centering}
\includegraphics[width=85mm]{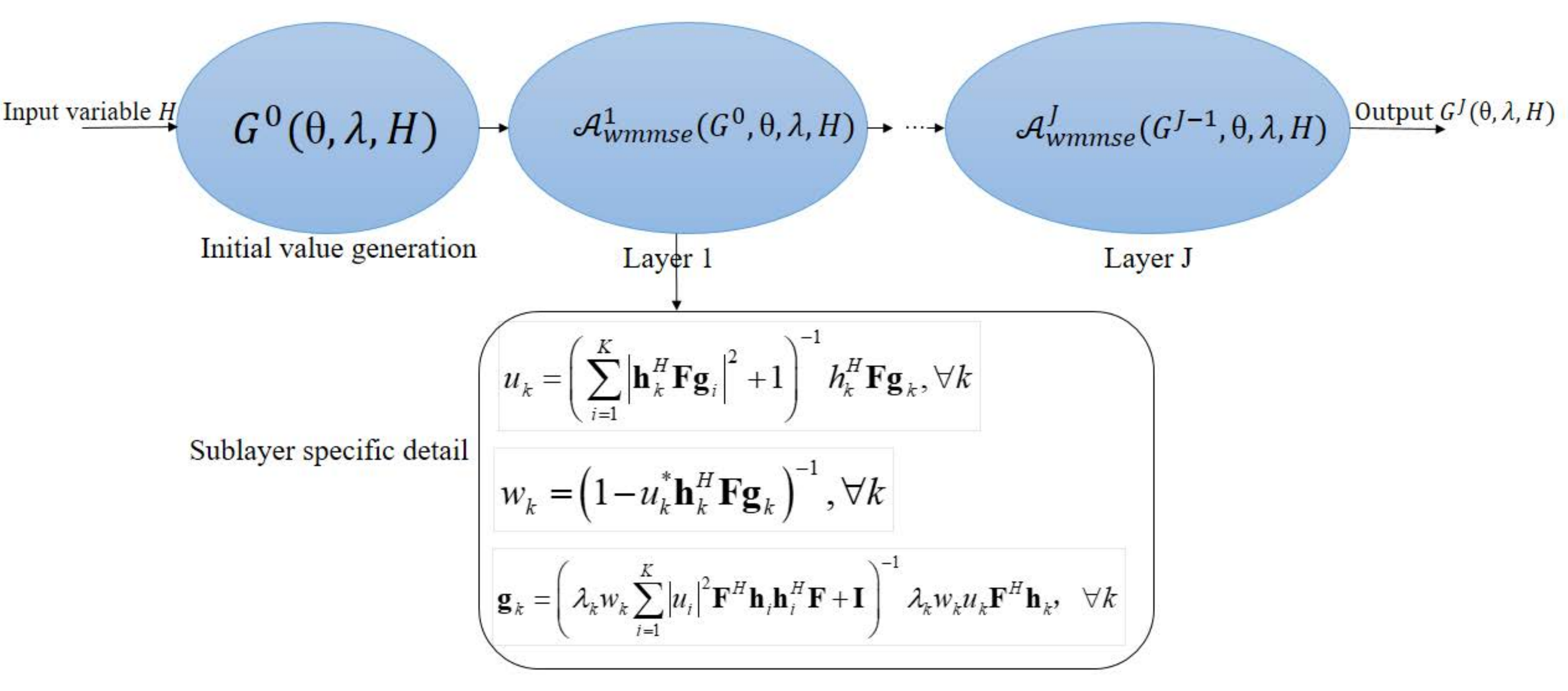}
\par\end{centering}
\caption{\label{fig:NN-WMMSE}The NN representation of the WMMSE short-term
sub-algorithm.}
\end{figure}

To be more specific, in Fig. \ref{fig:NN-WMMSE}, we illustrate the
NN representation of the WMMSE short-term sub-algorithm for application
example 2. In this case, the input is the channel state $\boldsymbol{H}$,
the output is the digital precoder $\boldsymbol{G}^{J}\left(\boldsymbol{\theta},\boldsymbol{\lambda},\boldsymbol{H}\right)$,
the $0$-th layer is used to generate the initial value $\boldsymbol{G}^{0}\left(\boldsymbol{\theta},\boldsymbol{\lambda},\boldsymbol{H}\right)$,
and the $j$-th layer $\mathcal{A}^{j}\left(\boldsymbol{G}^{j-1},\boldsymbol{\theta},\boldsymbol{\lambda},\boldsymbol{H}\right)$
is given by the update equations in (\ref{eq:Sw}) to (\ref{eq:sG}),
where $\boldsymbol{G}^{j-1}$ is an abbreviation for $\boldsymbol{G}^{j-1}\left(\boldsymbol{\theta},\boldsymbol{\lambda},\boldsymbol{H}\right)$. 

\begin{figure}
\begin{centering}
\includegraphics[width=85mm]{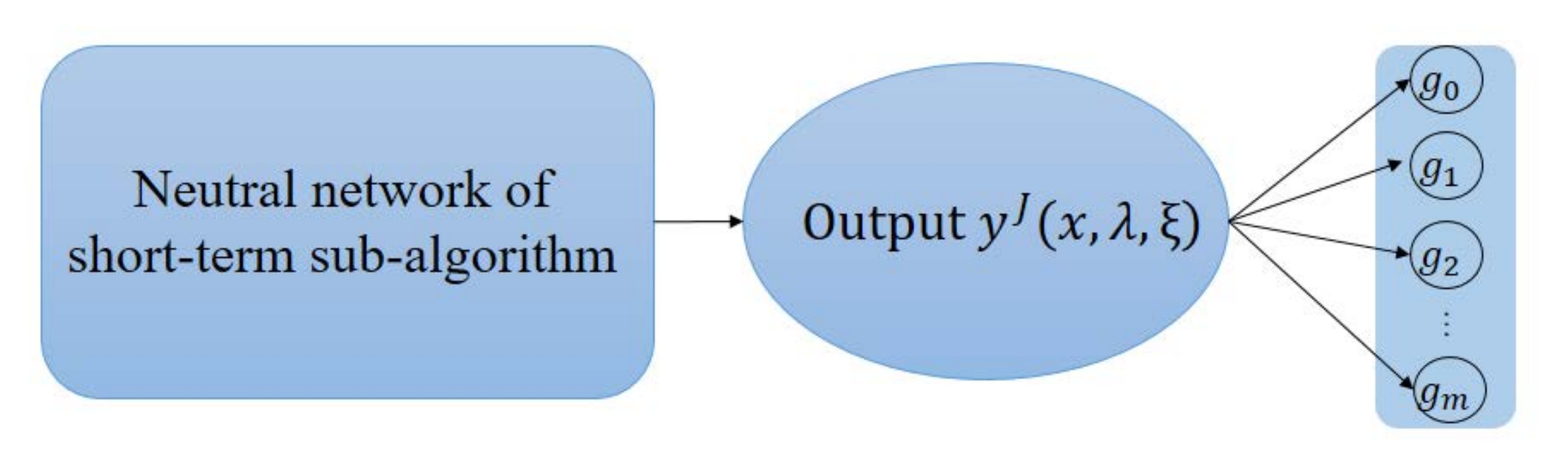}
\par\end{centering}
\caption{\label{fig:NN-g}NN to calculate the gradients of $g_{i}$'s.}
\end{figure}

With the NN representation of the short-term sub-algorithm, we can
use BP to calculate the gradients $\partial_{\boldsymbol{x}}\boldsymbol{y}^{J}\left(\boldsymbol{x},\boldsymbol{\lambda},\boldsymbol{\xi}\right)$
and $\partial_{\boldsymbol{\lambda}}\boldsymbol{y}^{J}\left(\boldsymbol{x},\boldsymbol{\lambda},\boldsymbol{\xi}\right)$
for the optimization of long-term variables. In fact, we can directly
use BP to calculate $\partial_{\boldsymbol{x}}\boldsymbol{y}^{J}\left(\boldsymbol{x},\boldsymbol{\lambda},\boldsymbol{\xi}\right)\partial_{\boldsymbol{y}}g_{i}\left(\boldsymbol{x},\boldsymbol{y}^{J}\left(\boldsymbol{x},\boldsymbol{\lambda},\boldsymbol{\xi}\right),\boldsymbol{\xi}\right),\forall i$
and $\partial_{\boldsymbol{\lambda}}\boldsymbol{y}^{J}\left(\boldsymbol{x},\boldsymbol{\lambda},\boldsymbol{\xi}\right)\partial_{\boldsymbol{y}}g_{i}\left(\boldsymbol{x},\boldsymbol{y}^{J}\left(\boldsymbol{x},\boldsymbol{\lambda},\boldsymbol{\xi}\right),\boldsymbol{\xi}\right),\forall i$
at $\left(\boldsymbol{x},\boldsymbol{\lambda}\right)=\left(\boldsymbol{x}^{t},\boldsymbol{\lambda}^{t}\right)$
by adding an additional layer to implement the function $g_{i}\left(\boldsymbol{x}^{t},\boldsymbol{y}^{J}\left(\boldsymbol{x},\boldsymbol{\lambda},\boldsymbol{\xi}\right),\boldsymbol{\xi}\right),\forall i$
on top of the NN representation of the short-term sub-algorithm, as
shown in Fig. \ref{fig:NN-g}. Note that the calculations for different
$i$ using BP share the same calculations over the NN of the short-term
sub-algorithm. With $\partial_{\boldsymbol{x}}\boldsymbol{y}^{J}\left(\boldsymbol{x},\boldsymbol{\lambda},\boldsymbol{\xi}\right)\partial_{\boldsymbol{y}}g_{i}\left(\boldsymbol{x},\boldsymbol{y}^{J}\left(\boldsymbol{x},\boldsymbol{\lambda},\boldsymbol{\xi}\right),\boldsymbol{\xi}\right),\forall i$
and $\partial_{\boldsymbol{x}}\boldsymbol{y}^{J}\left(\boldsymbol{x},\boldsymbol{\lambda},\boldsymbol{\xi}\right)\partial_{\boldsymbol{\lambda}}g_{i}\left(\boldsymbol{x},\boldsymbol{y}^{J}\left(\boldsymbol{x},\boldsymbol{\lambda},\boldsymbol{\xi}\right),\boldsymbol{\xi}\right),\forall i$
calculated by BP, the long-term variables $\boldsymbol{x},\boldsymbol{\lambda}$
can be optimized based on the statistics of the random state $\boldsymbol{\xi}$,
or a data set containing a large number of $N$ state samples $\boldsymbol{\xi}^{t},t=1,...,N$,
where the state samples can also be observed in an online manner.
Once the optimized long-term variables $\boldsymbol{x},\boldsymbol{\lambda}$
is obtained, the short-term variables $\boldsymbol{y}^{J}\left(\boldsymbol{x},\boldsymbol{\lambda},\boldsymbol{\xi}\right)$
for each random state $\boldsymbol{\xi}$ can be calculated in an
online manner using the NN representation.

\subsection{Discussions on Possible Extensions}

In this paper, we focus on the case when the short sub-algorithm is
designed based on the optimization theory and is then unrolled without
any modifications to facilitate the calculation of gradient. In practical
implementation, we may add some additional layers/parameters to the
NN unrolled from a short-term sub-algorithm (such as the GP, MM or
WMMSE), and represent the mapping $\boldsymbol{y}^{J}\left(\boldsymbol{x},\boldsymbol{\lambda},\boldsymbol{\xi}\right)$
using a more general NN as $\boldsymbol{y}^{J}\left(\boldsymbol{x},\boldsymbol{\lambda},\boldsymbol{\xi}\right)=\boldsymbol{\phi}\left(\boldsymbol{x},\boldsymbol{\lambda},\boldsymbol{\xi};\boldsymbol{\theta}\right)$,
where the vector $\boldsymbol{\theta}$ contains the additional parameters.
These additional layers/parameters can be optimized to potentially
speed up the convergence speed or even improve the performance. For
example, for a GP short-term sub-algorithm, we may change the update
equation in (\ref{eq:GP}) to
\[
\mathcal{A}^{j}\left(\boldsymbol{y}^{j-1}\right)=\mathbb{P}_{\mathcal{Y}\left(\boldsymbol{x},\boldsymbol{\lambda},\boldsymbol{\xi}\right)}\left[\boldsymbol{y}^{j-1}-\boldsymbol{\alpha}_{j}\circ\partial_{\boldsymbol{y}}g_{s}\left(\boldsymbol{x},\boldsymbol{\lambda},\boldsymbol{y}^{j-1},\boldsymbol{\xi}\right)\right],
\]
where the scalar step size $\alpha_{j}$ in the original GP update
equation (\ref{eq:GP}) is expanded to a vector step size $\boldsymbol{\alpha}_{j}$,
and $\circ$ denotes the Hadamard product. Since the vector step size
$\boldsymbol{\alpha}_{j}$ includes the scalar step size $\alpha_{j}$
as a special case, we may optimize the vector step sizes $\left\{ \boldsymbol{\alpha}_{j},j=1,...,J\right\} $
to achieve a better performance. 

After introducing the additional layers/parameters, both the parameter
$\boldsymbol{\theta}$ and the long-term variables $\boldsymbol{x},\boldsymbol{\lambda}$
can be optimized by solving the following modified long-term problem:
\begin{align}
\hat{\mathcal{P}}_{L}:\:\min_{\boldsymbol{x},\boldsymbol{\lambda},\boldsymbol{\theta}} & \hat{f}_{0}(\boldsymbol{x},\boldsymbol{\lambda},\boldsymbol{\theta})\triangleq\mathbb{E}\left[g_{0}\left(\boldsymbol{x},\boldsymbol{\phi}\left(\boldsymbol{x},\boldsymbol{\lambda},\boldsymbol{\xi};\boldsymbol{\theta}\right),\boldsymbol{\xi}\right)\right],\label{eq:mainP-1-2}\\
\text{s.t. } & \hat{f}_{i}(\boldsymbol{x},\boldsymbol{\lambda},\boldsymbol{\theta})\triangleq\mathbb{E}\left[g_{i}\left(\boldsymbol{x},\boldsymbol{\phi}\left(\boldsymbol{x},\boldsymbol{\lambda},\boldsymbol{\xi};\boldsymbol{\theta}\right),\boldsymbol{\xi}\right)\right]\leq0,\:\forall i.\nonumber 
\end{align}
Note that since the original short-term sub-algorithm $\boldsymbol{y}^{J}\left(\boldsymbol{x},\boldsymbol{\lambda},\boldsymbol{\xi}\right)$
is a special case of the NN $\boldsymbol{\phi}\left(\boldsymbol{x},\boldsymbol{\lambda},\boldsymbol{\xi};\boldsymbol{\theta}\right)$,
it is expected that a better objective value can be potentially achieved
by solving $\hat{\mathcal{P}}_{L}$ due to the extra freedom introduced
by the additional parameters $\boldsymbol{\theta}$. Problem $\hat{\mathcal{P}}_{L}$
can be solved using the same CSSCA method as in Algorithm 1, with
$\boldsymbol{\theta}$ as an additional long-term variable. Once the
optimized long-term variables $\boldsymbol{x},\boldsymbol{\lambda},\boldsymbol{\theta}$
is obtained, the short-term variables $\boldsymbol{y}^{J}\left(\boldsymbol{x},\boldsymbol{\lambda},\boldsymbol{\xi}\right)=\boldsymbol{\phi}\left(\boldsymbol{x},\boldsymbol{\lambda},\boldsymbol{\xi};\boldsymbol{\theta}\right)$
for each random state $\boldsymbol{\xi}$ can be calculated in an
online manner using the NN representation. Therefore, we may view
the Algorithm 1 as an unsupervised learning algorithm for the additional
parameters $\boldsymbol{\theta}$. It is also possible to pre-train
$\boldsymbol{\theta}$ using a supervised learning algorithm. Then
we may use the pre-trained $\boldsymbol{\theta}$ as an initial point
for the Algorithm 1. The pre-training helps to speed up the convergence
speed or even improve the performance. 

Finally, we give an example to briefly illustrate how to unfold the
WMMSE short-term sub-algorithm in Example 2, based on the deep-unfolding
framework proposed in \cite{Hu_TWC20_DeepunrollWMMSE}, where a general
form of iterative algorithm induced deep unfolding neural network
(IAIDNN) is developed in matrix form to better solve optimization
problems in communication systems. Specifically, the iterative WMMSE
algorithm is unfolded into a layer-wise structure with a series of
matrix multiplication and non-linear operations, as illustrated in
Fig. 2 in \cite{Hu_TWC20_DeepunrollWMMSE}. Compared to the original
WMMSE iterations, there are two major modifications. First, the element-wise
non-linear function and the first-order Taylor expansion structure
of the inverse matrix is applied to approximate the matrix inversion
operation in the WMMSE algorithm. Second, trainable parameters $\boldsymbol{\theta}$
are introduced in the non-linear function and the first-order Taylor
expansion structure of the inverse matrix, aiming at compensating
the performance loss caused by the approximation of matrix inversion.
Such modifications help to greatly speed up the convergence speed
and reduce the complexity of each iteration (layer) with little performance
loss. Then, the training process can be divided into the following
two stages. In the supervised learning stage, the expected distance
between the output of the IAIDNN $\boldsymbol{G}\left(\boldsymbol{x},\boldsymbol{\lambda},\boldsymbol{\xi}\right)$
and the labels produced by the iterative WMMSE algorithm $\boldsymbol{G}^{*}\left(\boldsymbol{x},\boldsymbol{\lambda},\boldsymbol{\xi}\right)$
under different long-term variables $\boldsymbol{x},\boldsymbol{\lambda}$
and channel realizations $\boldsymbol{\xi}$ is used as the loss function
to pre-train $\boldsymbol{\theta}$, with learning rate $\text{max}\left(0.1n^{-0.5},10^{-4}\right)$
and training batch size 100, where $n$ denotes the training iteration
number. Please refer to \cite{Hu_TWC20_DeepunrollWMMSE} for more
details. After applying the supervised learning, Algorithm 1 can be
used to optimize/train both the long-term variables $\boldsymbol{x},\boldsymbol{\lambda}$
and the parameter $\boldsymbol{\theta}$ in an unsupervised way.

\section{Applications\label{sec:Applications}}

In this section, we shall apply the proposed PDD-SSCA to solve the
two applications described in Section \ref{subsec:Examples-of-Problem}.
In both applications, we assume that the channel statistics is known
at the optimizer, and the \textit{coherence time of the channel statistics}
(within which the channel statistics is assumed to be constant) consists
of $T_{c}=1000$ channel realizations. The optimization algorithms
first calculate the optimal long-term variables at the beginning of
the channel statistics coherence time according to the known channel
statistics, and then solve the short-term subproblem to obtain the
short-term variables for each channel realization. The simulation
results are also obtained by averaging over $T_{c}=1000$ channel
realizations. The key algorithm parameters used in the simulations
are shown in Table \ref{tab:simulation-parameters}.

\begin{table}
\begin{centering}
\begin{tabular}{|c|c|c|}
\hline 
 & {\small{}Example 1} & {\small{}Example 2}\tabularnewline
\hline 
{\small{}Batch Size $B$} & {\small{}20} & {\small{}20}\tabularnewline
\hline 
{\small{}Step size $\rho^{t}$} & {\small{}$10/(10+n)^{0.9}$} & {\small{}$10/(10+n)^{0.9}$}\tabularnewline
\hline 
{\small{}Step size $\gamma^{t}$} & {\small{}$15\text{/}(15+n)$} & {\small{}$15\text{/}(15+n)$}\tabularnewline
\hline 
{\small{}Number of Layers $J$} & {\small{}N/A} & {\small{}$\ensuremath{5}$}\tabularnewline
\hline 
\end{tabular}
\par\end{centering}
\centering{}\caption{\textcolor{blue}{\label{tab:simulation-parameters}}Parameters of
the proposed PDD-SSCA algorithm used in the simulations.}
\end{table}

\subsection{Cognitive Multiple Access Channels}

In this example, the sample objective and constraint functions for
a fixed state $\boldsymbol{a},\boldsymbol{b}$ are given by
\begin{align}
g_{0}\left(\boldsymbol{\lambda},\Upsilon,\boldsymbol{a},\boldsymbol{b}\right) & =\log\left(1+\sum_{i=1}^{N}\frac{1}{N}\left(\frac{a_{i}}{b_{i}\Upsilon+\lambda_{i}}-1\right)^{+}\right),\nonumber \\
g_{i}\left(\boldsymbol{\lambda},\Upsilon,\boldsymbol{a},\boldsymbol{b}\right) & =\frac{1}{Na_{i}}\left(\frac{a_{i}}{b_{i}\Upsilon+\lambda_{i}}-1\right)^{+}-P_{i},\nonumber \\
g_{N+1}\left(\boldsymbol{\lambda},\Upsilon,\boldsymbol{a},\boldsymbol{b}\right) & =\sum_{i=1}^{N}\frac{b_{i}}{Na_{i}}\left(\frac{a_{i}}{b_{i}\Upsilon+\lambda_{i}}-1\right)^{+}-\Gamma,\label{eq:exam1gi}
\end{align}
for $i=1,...,N$, where $g_{i}\left(\boldsymbol{\lambda},\Upsilon,\boldsymbol{a},\boldsymbol{b}\right)$
is an abbreviation for $g_{i}\left(\boldsymbol{p}^{*}\left(\boldsymbol{\lambda},\Upsilon,\boldsymbol{a},\boldsymbol{b}\right),\boldsymbol{a},\boldsymbol{b}\right)$,
and $\boldsymbol{p}^{*}\left(\boldsymbol{\lambda},\Upsilon,\boldsymbol{a},\boldsymbol{b}\right)$
is the optimal short-term power allocation in (\ref{eq:Sexam1}).
To construct the surrogate function, we first calculate the gradients
$\nabla_{\boldsymbol{\lambda}}g_{i}\left(\boldsymbol{\lambda},\Upsilon,\boldsymbol{a},\boldsymbol{b}\right)$
and $\nabla_{\Upsilon}g_{i}\left(\boldsymbol{\lambda},\Upsilon,\boldsymbol{a},\boldsymbol{b}\right)$
for $i=0,1,....,N+1$. The calculation is straightforward from (\ref{eq:exam1gi})
and the details are omitted for conciseness. Note that strictly speaking,
$g_{i}$ is non-differentiable when $\frac{a_{i}}{b_{i}\Upsilon+\lambda_{i}}-1=0$.
However, since the probability of $\frac{a_{i}}{b_{i}\Upsilon+\lambda_{i}}-1=0$
is very small (zero for continuous state distribution), we can safely
ignore this point. With the gradients of $g_{i}$'s, we can construct
quadratic surrogate functions as
\begin{align}
\bar{f}_{i}^{t}\left(\boldsymbol{\lambda},\Upsilon\right) & =f_{i}^{t}+(\mathbf{f}_{\lambda,i}^{t})^{T}(\boldsymbol{\lambda}-\boldsymbol{\lambda}^{t})+(\mathbf{f}_{\Upsilon,i}^{t})^{T}(\Upsilon-\Upsilon^{t})\nonumber \\
 & +\tau_{i}\left(\left\Vert \boldsymbol{\lambda}-\boldsymbol{\lambda}^{t}\right\Vert ^{2}+\left(\Upsilon-\Upsilon^{t}\right)^{2}\right),\label{eq:SSF-exm1}
\end{align}
p'pwhere $\mathbf{f}_{\lambda,i}^{t}$ and $\mathbf{f}_{\Upsilon,i}^{t}$
can be calculated recursively as
\begin{align}
\mathbf{f}_{\Upsilon,i}^{t} & =\left(1-\rho^{t}\right)\mathbf{f}_{\Upsilon,i}^{t-1}+\rho^{t}\frac{1}{B}\sum_{j=1}^{B}\nabla_{\Upsilon}g_{i}\left(\boldsymbol{\lambda}^{t},\Upsilon^{t},\boldsymbol{a}_{j}^{t},\boldsymbol{b}_{j}^{t}\right),\nonumber \\
\mathbf{f}_{\lambda,i}^{t} & =\left(1-\rho^{t}\right)\mathbf{f}_{\lambda,i}^{t-1}+\rho^{t}\frac{1}{B}\sum_{j=1}^{B}\nabla_{\boldsymbol{\lambda}}g_{i}\left(\boldsymbol{\lambda}^{t},\Upsilon^{t},\boldsymbol{a}_{j}^{t},\boldsymbol{b}_{j}^{t}\right).\label{eq:fyx-exm1}
\end{align}
Note that the quadratic surrogate function in (\ref{eq:SSF-exm1})
is a special case of the structured surrogate function in (\ref{eq:SSF})
when $g_{i}^{c}\left(\boldsymbol{x},\boldsymbol{y},\boldsymbol{\xi}\right)$
is set to be 0. In this case, both the objective update (\ref{eq:Pitert})
and the feasible update (\ref{eq:Pitert-1}) are simple quadratic
optimization problems, which can be solved efficiently using the existing
optimization solvers such as CVX \cite{cvx}.

We compare the proposed PDD-SSCA with the following baseline algorithms.
\begin{itemize}
\item \textbf{Baseline 1 (DL)}: This is the DL algorithm in \cite{Quek_JSAC2019_DLDist}. 
\item \textbf{Baseline 2 (Dual Ellipsoid)}: This is the dual decomposition
algorithm in \cite{Zhang_TIT2009_CMAC}. The dual problem is solved
using the Ellipsoid method \cite{Bland_Oper81_ellipsoid}. Note that
the dual decomposition is applicable to Example 1 because this problem
is convex. In each iteration of the Ellipsoid method, the expectation
$\mathbb{E}\left[g_{i}\left(\boldsymbol{\lambda},\Upsilon,\boldsymbol{a},\boldsymbol{b}\right)\right],\forall i$
is required to calculate the subgradient, which is obtained by sample
averaging over $T_{c}=1000$ channel realizations.
\item \textbf{Baseline 3 (Short-term Constraint)}: Example Problem 1 is
optimally solved with the short-term constraints $p_{i}\left(\boldsymbol{a},\boldsymbol{b}\right)\leq P_{i}$
and $\sum_{i=1}^{N}b_{i}p_{i}\left(\boldsymbol{a},\boldsymbol{b}\right)\leq\Gamma$,
instead of $\mathbb{E}\left[p_{i}\left(\boldsymbol{a},\boldsymbol{b}\right)\right]\leq P_{i}$
and $\mathbb{E}\left[\sum_{i=1}^{N}b_{i}p_{i}\left(\boldsymbol{a},\boldsymbol{b}\right)\right]\leq\Gamma$.
\end{itemize}

\begin{figure}
\begin{centering}
\includegraphics[width=85mm]{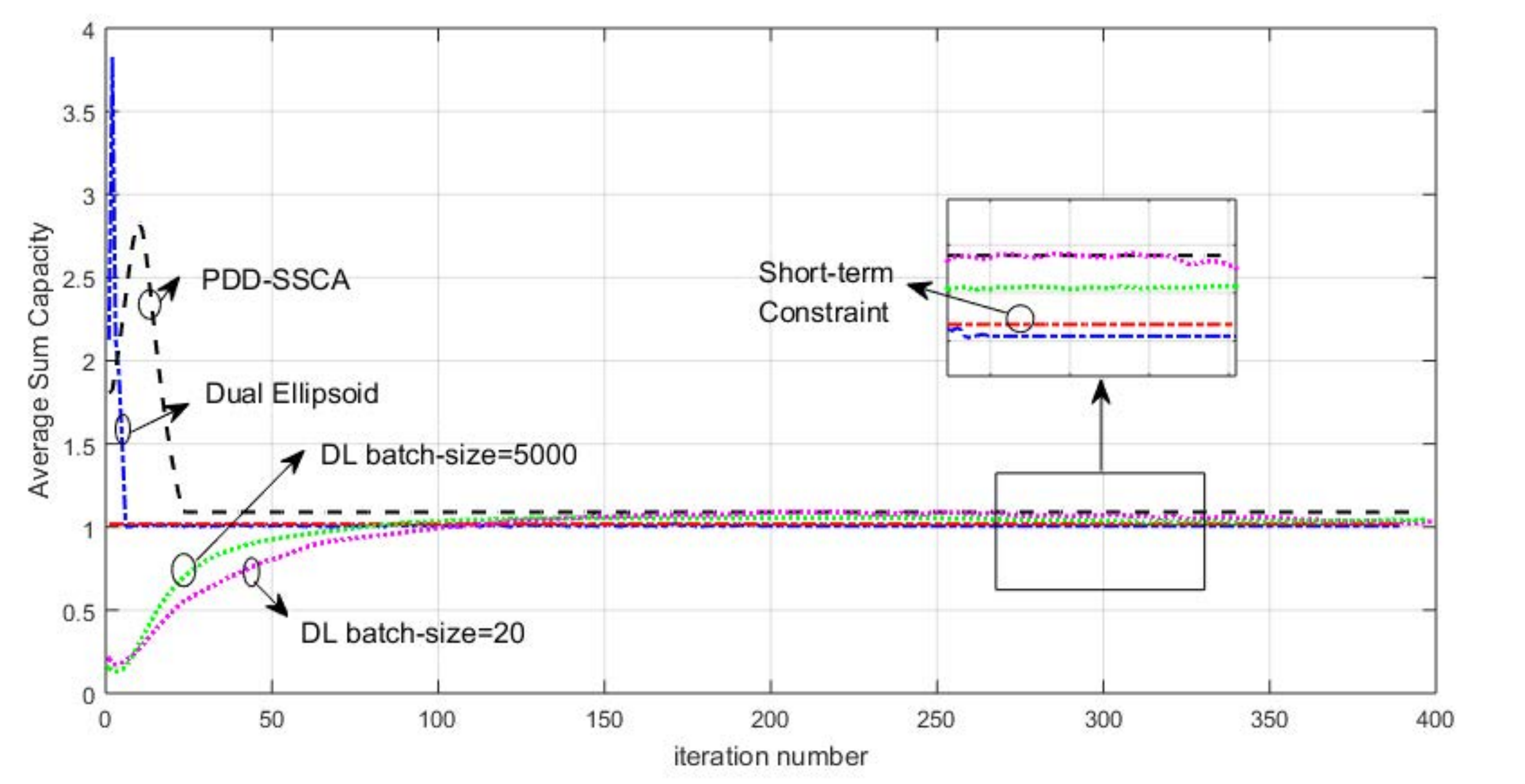}
\par\end{centering}
\caption{\label{fig:exam1f1}Average sum capacity versus iteration number.}
\end{figure}

\begin{figure}
\begin{centering}
\includegraphics[width=85mm]{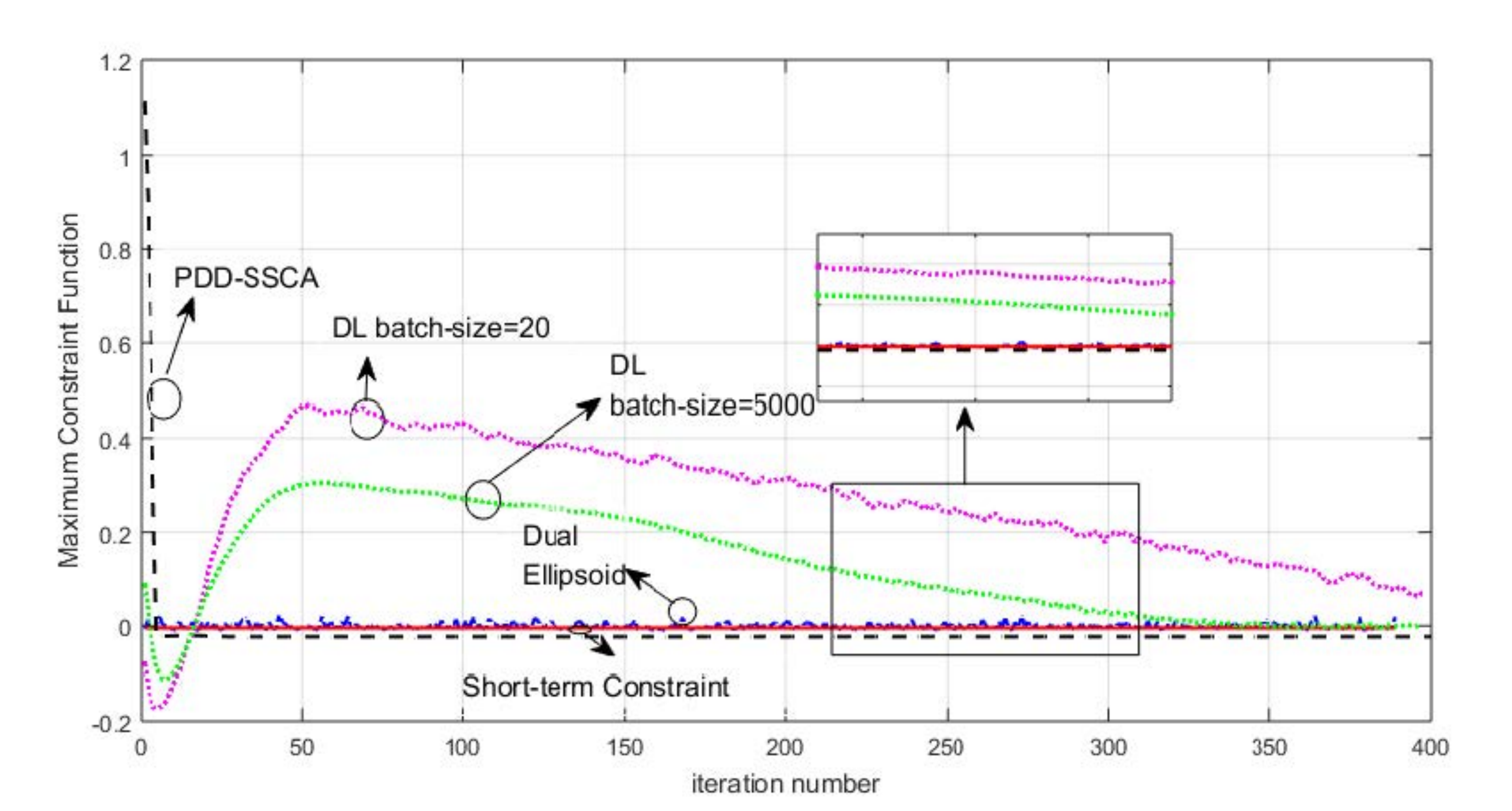}
\par\end{centering}
\caption{\label{fig:exam1f2}Maximum constraint function versus iteration number.}
\end{figure}

\begin{table}
\begin{centering}
{\footnotesize{}}%
\begin{tabular}{|l|l|l|}
\hline 
 & {\small{}Complexity order per iteration} & {\small{}Mem. cost}\tabularnewline
\hline 
{\small{}Proposed} & {\small{}$O(I_{1}(NlogN+B_{1}N)+T_{c}N)$} & {\small{}$O(B_{1}N)$ }\tabularnewline
{\small{}DL} & {\small{}$O\left(I_{2}B_{2}100N^{2}L+T_{c}N\right)$} & {\small{}$O(B_{2}N+$ }\tabularnewline
 &  & {\small{}$100N^{2}L)$}\tabularnewline
{\small{}Dual Ellipsoid} & {\small{}$O\left(I_{3}(N^{2}+T_{c}N\right)+T_{c}N)$} & {\small{}$O(NT_{c})$ }\tabularnewline
{\small{}Baseline 3} & {\small{}$O\left(T_{c}N^{2}\right)$} & {\small{}$O\left(N\right)$ }\tabularnewline
\hline 
\end{tabular}{\footnotesize\par}
\par\end{centering}
{\small{}\caption{\label{tab:Cputime}Comparison of the complexity order and memory
cost for different algorithms in Example 1.}
}{\small\par}
\end{table}

Following the same simulation configuration as in \cite{Quek_JSAC2019_DLDist},
we set $N=2$, $\Gamma=0.5$, $P_{i}=P=5\text{ dB},\forall i$. In
Figs. \ref{fig:exam1f1} and \ref{fig:exam1f2}, we plot the objective
function (average sum capacity) and maximum constraint function versus
the iteration number, respectively. The per-iteration computation
complexity and memory cost are compared in Table \ref{tab:Cputime},
where the $B_{i},i=1,2$ and $I_{i},i=1,2,3$ account for the mini-batch
size and the number of iterations for optimizing the long-term variables,
respectively. The complexity of calculating the estimated long-term
gradient is $O\left(B_{1}N\right)$ for mini-batch size $B_{1}$.
The complexity of solving the objective update problem in (\ref{eq:Pitert})
or the feasibility update problem in (\ref{eq:Pitert-1}) using the
interior-point method is given by $O\left(N\left(log\left(N\right)\right)\right)$\footnote{This is because the objective/feasibility update problem is a convex
quadratic programming problem with a diagonal Hessian matrix.}. Moreover, the complexity of calculating the optimal short-term solution
is $O\left(N\right)$ per channel realization. Therefore, the per-iteration
complexity order of the proposed algorithm is as given in Table \ref{tab:Cputime}.
The memory cost is proportional to the amount of parameters to be
optimized plus the product of the batch size and the dimension of
the state variable, i.e., $O(B_{1}N)$ Bytes. The complexity order
and memory cost of the baseline algorithms can be analyzed similarly.
PDD-SSCA converges to the optimal average sum capacity with all average
transmit power constraints and interference threshold constraint satisfied
with high accuracy. For the same batch size 20, the number of iterations
required to achieve a good convergence accuracy in the proposed PDD-SSCA
is much less than that in the DL (25 versus more than 400). Moreover,
the complexity/memory cost of the PDD-SSCA is also much less than
that of the DL since it exploits the structure of the problem and
has much less parameters to be optimized (3 versus 2080). When the
batch size of the DL is increased to 5000, the required number of
iterations of DL can be reduced (but still much larger than 25), at
the cost of dramatically increasing the complexity and memory cost.
Both the convergence speed and performance of the PDD-SSCA are similar
to that of the Dual Ellipsoid method. However, the complexity and
memory cost of the PDD-SSCA are much lower. Therefore, the proposed
PDD-SSCA is much more efficient than the DL and Dual Ellipsoid methods.
Finally, the complexity/memory cost of the PDD-SSCA is similar to/slightly
higher than the ``Short-term Constraint'' baseline, but the performance
is much better. 

\subsection{Power Minimization for Two-timescale Hybrid Beamforming}

In this example, the short-term WMMSE algorithm can be implemented
using a simple NN as shown in Fig. \ref{fig:NN-WMMSE}. Moreover,
the gradients of the sample objective and constraint functions $\nabla_{\boldsymbol{\theta}}g_{i}\left(\boldsymbol{\theta},\boldsymbol{G}^{J}\left(\boldsymbol{\theta},\boldsymbol{\lambda},\boldsymbol{H}\right),\boldsymbol{H}\right)$
and $\nabla_{\boldsymbol{\lambda}}g_{i}\left(\boldsymbol{\theta},\boldsymbol{G}^{J}\left(\boldsymbol{\theta},\boldsymbol{\lambda},\boldsymbol{H}\right),\boldsymbol{H}\right)$
for a fixed state $\boldsymbol{H}$ can be calculated using the BP
method based on the NN in Fig. \ref{fig:NN-g}. Then, we can construct
quadratic surrogate functions as
\begin{align}
\bar{f}_{i}^{t}\left(\boldsymbol{\theta},\boldsymbol{\lambda}\right) & =f_{i}^{t}+(\mathbf{f}_{\boldsymbol{\theta},i}^{t})^{T}(\boldsymbol{\theta}-\boldsymbol{\theta}^{t})+(\mathbf{f}_{\lambda,i}^{t})^{T}(\boldsymbol{\lambda}-\boldsymbol{\lambda}^{t})\nonumber \\
 & +\tau_{i}\left(\left\Vert \boldsymbol{\theta}-\boldsymbol{\theta}^{t}\right\Vert ^{2}+\left\Vert \boldsymbol{\lambda}-\boldsymbol{\lambda}^{t}\right\Vert ^{2}\right),\label{eq:SSF-exm2}
\end{align}
where $\mathbf{f}_{\boldsymbol{\theta},i}^{t}$ and $\mathbf{f}_{\lambda,i}^{t}$
can be calculated recursively as 
\begin{align}
\mathbf{f}_{\boldsymbol{\theta},i}^{t} & =\left(1-\rho^{t}\right)\mathbf{f}_{\boldsymbol{\theta},i}^{t-1}+\rho^{t}\frac{1}{B}\sum_{j=1}^{B}\nabla_{\boldsymbol{\theta}}g_{i}\left(\boldsymbol{\theta}^{t},\boldsymbol{\lambda}^{t},\boldsymbol{H}_{j}^{t}\right),\nonumber \\
\mathbf{f}_{\lambda,i}^{t} & =\left(1-\rho^{t}\right)\mathbf{f}_{\lambda,i}^{t-1}+\rho^{t}\frac{1}{B}\sum_{j=1}^{B}\nabla_{\boldsymbol{\lambda}}g_{i}\left(\boldsymbol{\theta}^{t},\boldsymbol{\lambda}^{t},\boldsymbol{H}_{j}^{t}\right),\label{eq:fyx-exm2}
\end{align}
where $g_{i}\left(\boldsymbol{\theta}^{t},\boldsymbol{\lambda}^{t},\boldsymbol{H}_{j}^{t}\right)$
is an abbreviation for $g_{i}\left(\boldsymbol{\theta}^{t},\boldsymbol{G}^{J}\left(\boldsymbol{\theta}^{t},\boldsymbol{\lambda}^{t},\boldsymbol{H}_{j}^{t}\right),\boldsymbol{H}_{j}^{t}\right)$.
Again, both the objective update (\ref{eq:Pitert}) and the feasible
update (\ref{eq:Pitert-1}) are simple quadratic optimization problems,
which can be solved efficiently. 

We compare the proposed PDD-SSCA with the following baseline algorithms.
\begin{itemize}
\item \textbf{Baseline 1 (SSCA-THP)}: The SSCA-THP in \cite{Liu_JSTSP2018_SSCATHP}
is a single-stage optimization algorithm which only optimizes the
RF precoder $\boldsymbol{\theta}$ with the baseband precoder fixed
as the regularized zero-forcing (RZF) precoder \cite{Peel_TOC05_RCI}.
\item \textbf{Baseline 2 (Prima-Dual)}: This is a primal-dual method based
algorithm, in which the primal-dual method in \cite{Quek_JSAC2019_DLDist}
is employed to solve the long-term problem in (\ref{eq:mainP-1})
and the rest is similar to the proposed algorithm.
\item \textbf{Baseline 3 (SAA-PDD-SCA)}: Sample average approximation (SAA)
is a common method to solve a stochastic optimization problem \cite{SPlecbook},
at the cost of high computation and memory overhead. After applying
the SAA on the constraint functions using $200$ channel samples,
the problem becomes a deterministic optimization problem with coupled
non-convex constraints. We apply the proposed primal-dual decomposition
to decompose this deterministic problem into one master problem (optimization
of RF precoding and long-term Lagrange multipliers) and 200 subproblems
(optimization of digital precoding). Then in each iteration, the master
problem is solved using the deterministic SCA method in \cite{Meisam_thesis14_BSUM}
and the digital precoding optimization subproblems are solved using
the WMMSE method.
\item \textbf{Baseline 4 (SLNR-max)}: This is the SLNR-max algorithm in
\cite{Park_TSP17_THP}. The RF precoder is optimized by maximizing
the signal-to-leakage plus noise ratio (SLNR) \cite{Park_TSP17_THP}.
Then the digital precoder for each channel realization is optimized
using the WMMSE method by solving a power minimization problem subject
to the rate constraint.
\end{itemize}

\begin{figure}
\begin{centering}
\includegraphics[width=85mm]{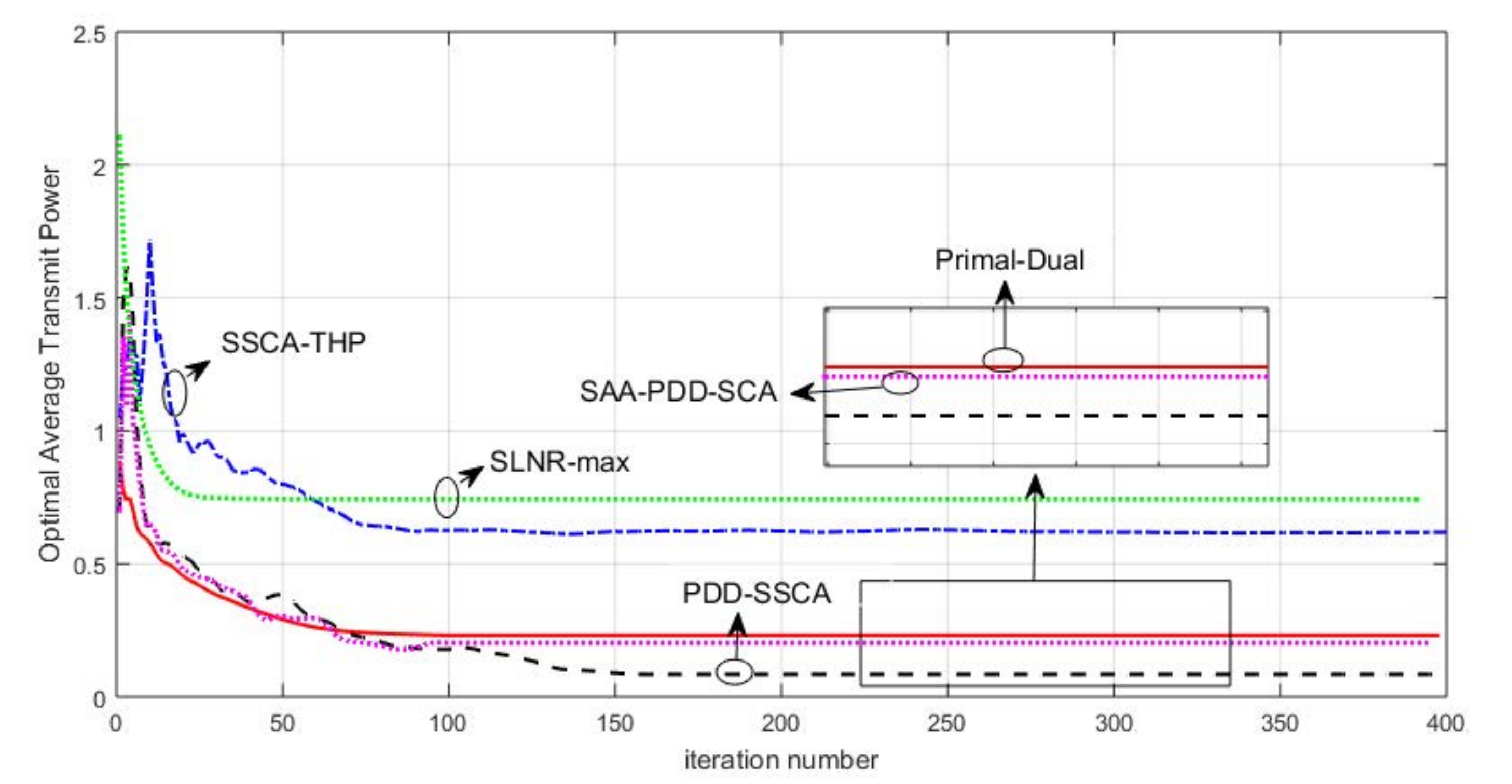}
\par\end{centering}
\caption{\label{fig:exam2f1}Average transmit power versus iteration number.}
\end{figure}

\begin{figure}
\begin{centering}
\includegraphics[width=85mm]{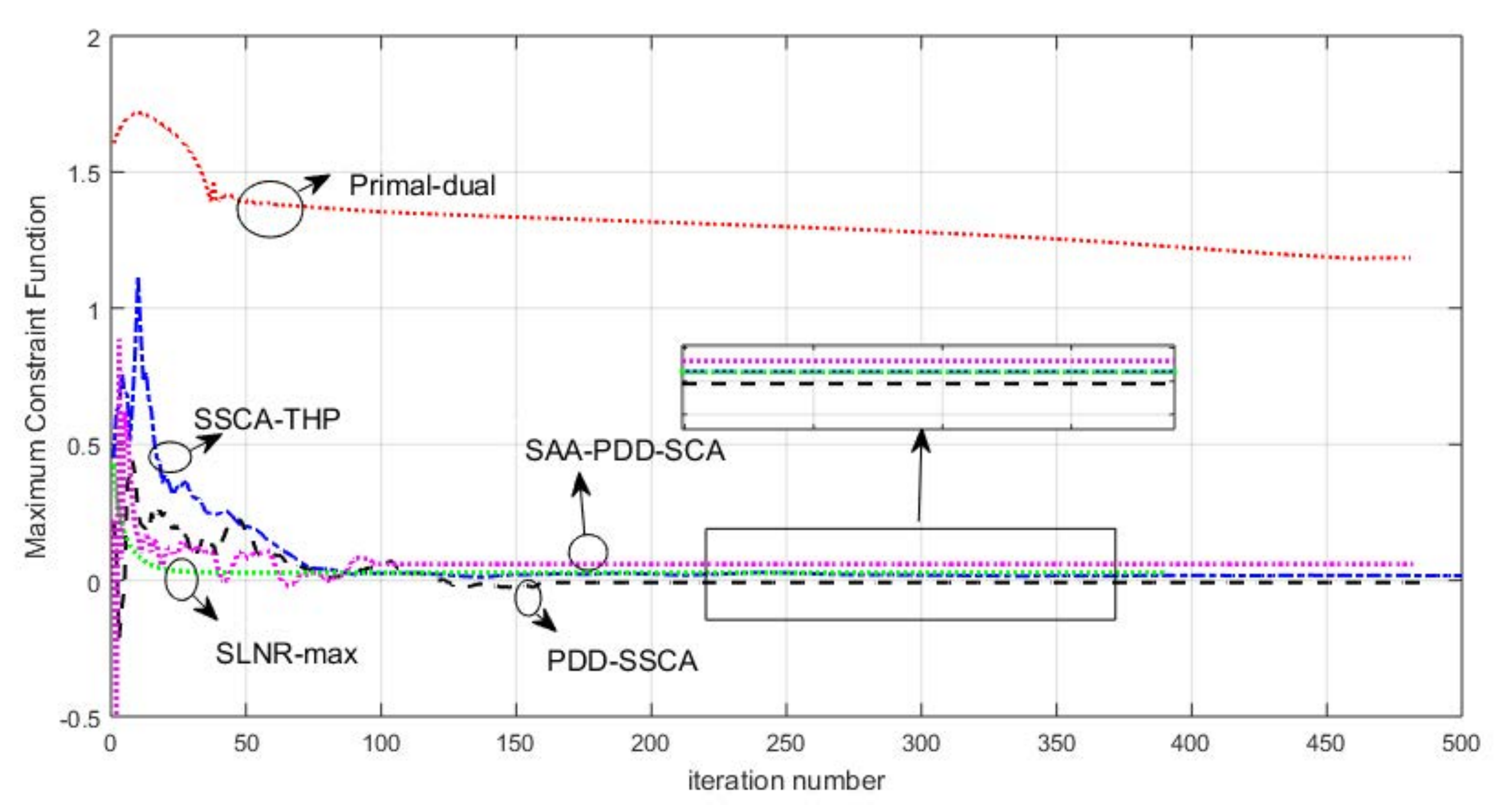}
\par\end{centering}
\caption{\label{fig:exam2f2} Maximum constraint function versus iteration
number.}
\end{figure}

\begin{table*}
\begin{centering}
{\footnotesize{}}%
\begin{tabular}{|l|l|l|}
\hline 
 & {\small{}Complexity order per iteration} & {\small{}Memory cost}\tabularnewline
\hline 
{\small{}Proposed} & {\small{}$O(I_{1}(MSlogK+B_{1}I_{W}S^{3})+T_{c}I_{W}S^{3})$} & {\small{}$O(MS+B_{1}MK)$}\tabularnewline
{\small{}SSCA-THP} & {\small{}$O\left(I_{2}\left(MSlog\left(K\right)+B_{1}S^{3}\right)+T_{c}S^{3}\right)$} & {\small{}$O(MS+B_{1}MK)$ }\tabularnewline
{\small{}Prima-Dual} & {\small{}$O\left(I_{3}\left(MS+B_{1}I_{W}S^{3}\right)+T_{c}I_{W}S^{3}\right)$} & {\small{}$O(MS+B_{1}MK)$}\tabularnewline
{\small{}SAA-PDD-SCA} & {\small{}$O(I_{4}(MSlogK+B_{2}I_{W}S^{3})+T_{c}I_{W}S^{3})$} & {\small{}$O(MS+B_{2}MK)$ }\tabularnewline
{\small{}SLNR-max} & {\small{}$O\left(\ensuremath{I_{5}M^{3}}+T_{c}I_{W}S^{3}\right)$} & {\small{}$O\left(KM^{2}\right)$ }\tabularnewline
\hline 
\end{tabular}{\footnotesize\par}
\par\end{centering}
{\small{}\caption{\label{tab:Cputime-1}Comparison of the complexity order and memory
cost for different algorithms in Example 2.}
}{\small\par}
\end{table*}

In the simulations, we adopt the same geometry-based channel model
as in \cite{Park_TSP17_THP,Liu_JSTSP2018_SSCATHP}. There are $M=64$
antennas and $S=4$ RF chains at the BS, serving $K=4$ users. In
Figs. \ref{fig:exam2f1} and \ref{fig:exam2f2}, we plot the objective
function (average transmit power) and maximum constraint function
(target average rate minus achieved average rate) versus the iteration,
respectively. The per-iteration computation complexity and memory
cost are compared in Table \ref{tab:Cputime-1}, where the $B_{i}$
and $I_{i}$ account for the mini-batch size and the number of iterations
for optimizing the long-term variables, $I_{W}$ is the number iterations
of WMMSE algorithm. We set $B_{1}=20$ and $B_{2}=200$ in the simulations.
The complexity order of the short-term WMMSE algorithm is given by
$O\left(I_{W}S^{3}\right)$. The complexity of calculating the estimated
long-term gradient is $O\left(B_{1}I_{W}S^{3}\right)$ for mini-batch
size $B_{1}$. The complexity of solving the objective/feasibility
update problem using the interior-point method is given by $O\left(MSlog\left(K\right)\right)$.
Moreover, the complexity of calculating the optimal short-term solution
is $O\left(I_{W}S^{3}\right)$ per channel realization. Therefore,
the per-iteration complexity order of the proposed algorithm is as
given in Table \ref{tab:Cputime-1}. The memory cost is proportional
to the amount of parameters to be optimized plus the product of the
batch size and the dimension of the state variable, i.e., $O(MS+B_{1}MK)$
Bytes. The complexity order and memory cost of the baseline algorithms
can be analyzed similarly. The number of iteration required to achieve
a good convergence accuracy in the proposed PDD-SSCA is similar to
that in the SSCA-THP. However, the PDD-SSCA converges to a much lower
average transmit power with all target average rates satisfied with
high accuracy. The PDD-SSCA can achieve a large performance gain over
the SSCA-THP. On the other hand, the primal-dual method based algorithm
cannot converge to a feasible solution that satisfies all the average
rate constraints, due to the highly non-convex nature of the constraints.
The convergence speed of the PDD-SSCA is similar to that of the SAA-PDD-SCA
method. However, the performance of the PDD-SSCA is better because
using $B_{2}=200$ samples is not sufficient to obtain a good sample
average approximation, and the complexity/memory cost of the PDD-SSCA
is also much lower. Finally, the complexity/memory cost of the PDD-SSCA
is lower than the SLNR-max baseline, and the performance is much better. 

\subsection{Conclusions\label{sec:Conclusion}}

We propose a PDD-SSCA algorithmic framework to solve a class of two-stage
stochastic optimization problems, in which the long-term and short-term
variables are tightly coupled in non-convex stochastic constraints.
The PDD-SSCA is designed based on a novel two-stage primal-dual decomposition
method established in this paper, which shows that the tightly coupled
two-stage problem can be decomposed into a long-term problem and a
family of short-term subproblems. At each iteration, PDD-SSCA first
runs a short-term sub-algorithm to find stationary points of the short-term
subproblems associated with a mini-batch of the state samples. Then
it constructs a convex surrogate for the long-term problem based on
the deep unrolling of the short-term sub-algorithm and the back propogation
method. Finally, the optimal solution of the convex surrogate problem
is solved to generate the next iterate. We show that under some technical
conditions, PDD-SSCA converges to a KKT solution of the original two-stage
problem almost surely. The effectiveness of the proposed PDD-SSCA
method is verified using two important application examples.

\appendix

\subsubsection{Proof of Theorem \ref{thm:PDidea} \label{subsec:Proof-of-Theorem-PDDidea}}

Using the condition $\delta_{i}\in\left(0,G_{i}^{\textrm{max}}\left(\boldsymbol{x}^{\circ}\right)-G_{i}^{\textrm{min}}\left(\boldsymbol{x}^{\circ}\right)\right),\forall i$
and following similar analysis as that in the proof of Proposition
2 of \cite{Quek_JSAC2019_DLDist}, it can be shown that for fixed
$\boldsymbol{x}^{\circ}$, Problem $\mathcal{P}$ fulfils the time-sharing
condition \cite{Yu_06_TCOM_Dual_nonconvex_opt} and the strong duality
holds for $\mathcal{P}$ with fixed $\boldsymbol{x}^{\circ}$. Therefore,
there must exist Lagrange multipliers $\boldsymbol{\lambda}^{\circ}$,
such that $\Theta^{\circ}$ is an optimal solution of the following
problem: 
\begin{align}
\min_{\Theta} & \mathbb{E}\left[g_{0}\left(\boldsymbol{x}^{\circ},\boldsymbol{y}\left(\boldsymbol{\xi}\right),\boldsymbol{\xi}\right)\right]\label{eq:mainPS-1}\\
 & +\sum_{i}\lambda_{i}^{\circ}\mathbb{E}\left[g_{i}\left(\boldsymbol{x}^{\circ},\boldsymbol{y}\left(\boldsymbol{\xi}\right),\boldsymbol{\xi}\right)\right],\nonumber \\
\text{s.t. } & h_{j}\left(\boldsymbol{y}\left(\boldsymbol{\xi}\right),\boldsymbol{\xi}\right)\leq0,\:j=1,...,n,\forall\boldsymbol{\xi}.\nonumber 
\end{align}
Moreover, 
\begin{align}
\lambda_{i}^{\circ}\mathbb{E}\left[g_{i}\left(\boldsymbol{x}^{\circ},\boldsymbol{y}^{\circ}\left(\boldsymbol{\xi}\right),\boldsymbol{\xi}\right)\right] & =0,\forall i,\nonumber \\
\mathbb{E}\left[g_{i}\left(\boldsymbol{x}^{\circ},\boldsymbol{y}^{\circ}\left(\boldsymbol{\xi}\right),\boldsymbol{\xi}\right)\right] & \leq0,\forall i.\label{eq:T1}
\end{align}
By definition, $\left\{ \boldsymbol{y}^{\star}\left(\boldsymbol{x}^{\circ},\boldsymbol{\lambda}^{\circ},\boldsymbol{\xi}\right),\forall\boldsymbol{\xi}\right\} $
is also an optimal solution of (\ref{eq:mainPS-1}) and satisfies
\begin{align}
\lambda_{i}^{\circ}\mathbb{E}\left[g_{i}\left(\boldsymbol{x}^{\circ},\boldsymbol{y}^{\star}\left(\boldsymbol{x}^{\circ},\boldsymbol{\lambda}^{\circ},\boldsymbol{\xi}\right),\boldsymbol{\xi}\right)\right] & =0,\forall i,\nonumber \\
\mathbb{E}\left[g_{i}\left(\boldsymbol{x}^{\circ},\boldsymbol{y}^{\star}\left(\boldsymbol{x}^{\circ},\boldsymbol{\lambda}^{\circ},\boldsymbol{\xi}\right),\boldsymbol{\xi}\right)\right] & \leq0,\forall i.\label{eq:T2}
\end{align}
From (\ref{eq:T1}), (\ref{eq:T2}) and the fact that both $\Theta^{\circ}$
and $\left\{ \boldsymbol{y}^{\star}\left(\boldsymbol{x}^{\circ},\boldsymbol{\lambda}^{\circ},\boldsymbol{\xi}\right),\forall\boldsymbol{\xi}\right\} $
are optimal solutions of (\ref{eq:mainPS-1}), we have
\begin{align*}
f_{0}^{\star}(\boldsymbol{x}^{\circ},\boldsymbol{\lambda}^{\circ}) & =\mathbb{E}\left[g_{0}\left(\boldsymbol{x}^{\circ},\boldsymbol{y}^{\star}\left(\boldsymbol{x}^{\circ},\boldsymbol{\lambda}^{\circ},\boldsymbol{\xi}\right),\boldsymbol{\xi}\right)\right].\\
 & =\mathbb{E}\left[g_{0}\left(\boldsymbol{x}^{\circ},\boldsymbol{y}^{\circ}\left(\boldsymbol{\xi}\right),\boldsymbol{\xi}\right)\right]=f_{0}(\boldsymbol{x}^{\circ},\Theta^{\circ}),\\
f_{i}^{\star}(\boldsymbol{x}^{\circ},\boldsymbol{\lambda}^{\circ}) & =\mathbb{E}\left[g_{i}\left(\boldsymbol{x}^{\circ},\boldsymbol{y}^{\star}\left(\boldsymbol{x}^{\circ},\boldsymbol{\lambda}^{\circ},\boldsymbol{\xi}\right),\boldsymbol{\xi}\right)\right]\leq0,\forall i.
\end{align*}
Therefore, $\boldsymbol{x}^{\circ},\boldsymbol{\lambda}^{\circ}$
is a feasible solution of $\mathcal{P}_{L}$, and thus $f_{0}^{\star}(\boldsymbol{x}^{\star},\boldsymbol{\lambda}^{\star})\leq f_{0}^{\star}(\boldsymbol{x}^{\circ},\boldsymbol{\lambda}^{\circ})=f_{0}(\boldsymbol{x}^{\circ},\Theta^{\circ})$,
from which it follows that $\left(\boldsymbol{x}^{\star},\Theta^{\star}\right)$
is also the optimal solution of $\mathcal{P}$.

\subsubsection{Proof of Theorem \ref{thm:PDprac} \label{subsec:Proof-of-Theorem-PDDprac}}

We first prove a useful Lemma.
\begin{lem}
\label{lem:keylemma1}Let $\mathcal{I}_{A}^{S}\left(\boldsymbol{\xi}\right)=\left\{ j:\nu_{j}\left(\boldsymbol{x}^{*},\boldsymbol{\lambda}^{*},\boldsymbol{\xi}\right)>0\right\} $.
We must have $\lim_{J\rightarrow\infty}\nabla_{\boldsymbol{\lambda}}h_{j}\left(\boldsymbol{y}^{J*},\boldsymbol{\xi}\right)=\boldsymbol{0},\forall j\in\mathcal{I}_{A}^{S}\left(\boldsymbol{\xi}\right)$,
where $\boldsymbol{y}^{J*}$ is an abbreivation for $\boldsymbol{y}^{J}\left(\boldsymbol{x}^{*},\boldsymbol{\lambda}^{*},\boldsymbol{\xi}\right)$.
\end{lem}
\begin{IEEEproof}
From the first KKT condition in (\ref{eq:STsconeJ}), and the LIRC,
it can be shown that $\nu_{j}\left(\boldsymbol{x},\boldsymbol{\lambda},\boldsymbol{\xi}\right)$
is differentiable w.r.t. $\boldsymbol{\lambda}$ and $\partial_{\boldsymbol{\lambda}}\nu_{j}\left(\boldsymbol{x}^{*},\boldsymbol{\lambda}^{*},\boldsymbol{\xi}\right)$
is bounded. From the third KKT condition in (\ref{eq:STsconeJ}),
we have 
\begin{align}
\nabla_{\boldsymbol{\lambda}}h_{j}(\boldsymbol{y}^{J*},\boldsymbol{\xi}) & =\frac{\partial_{\boldsymbol{\lambda}}e_{3,j}^{J}(\boldsymbol{x}^{*},\boldsymbol{\lambda}^{*},\boldsymbol{\xi})}{\nu_{j}\left(\boldsymbol{x}^{*},\boldsymbol{\lambda}^{*},\boldsymbol{\xi}\right)}\nonumber \\
 & -\frac{e_{3,j}^{J}(\boldsymbol{x}^{*},\boldsymbol{\lambda}^{*},\boldsymbol{\xi})\partial_{\boldsymbol{\lambda}}\nu_{j}(\boldsymbol{x}^{*},\boldsymbol{\lambda}^{*},\boldsymbol{\xi})}{\nu_{j}^{2}\left(\boldsymbol{x}^{*},\boldsymbol{\lambda}^{*},\boldsymbol{\xi}\right)},\label{eq:L1eq}
\end{align}
$\forall j\in\mathcal{I}_{A}^{S}\left(\boldsymbol{\xi}\right)$. Then
Lemma \ref{lem:keylemma1} follows from (\ref{eq:L1eq}) and the fact
that $\nu_{j}\left(\boldsymbol{x}^{*},\boldsymbol{\lambda}^{*},\boldsymbol{\xi}\right)>0,\forall j\in\mathcal{I}_{A}^{S}\left(\boldsymbol{\xi}\right)$,
$\partial_{\boldsymbol{\lambda}}\nu_{j}\left(\boldsymbol{x}^{*},\boldsymbol{\lambda}^{*},\boldsymbol{\xi}\right)$
is bounded and $\lim_{J\rightarrow\infty}e_{3,j}^{J}\left(\boldsymbol{x}^{*},\boldsymbol{\lambda}^{*},\boldsymbol{\xi}\right)=0$,
$\lim_{J\rightarrow\infty}\partial_{\boldsymbol{\lambda}}e_{3,j}^{J}\left(\boldsymbol{x}^{*},\boldsymbol{\lambda}^{*},\boldsymbol{\xi}\right)=\boldsymbol{0}$.
\end{IEEEproof}

According to the chain rule,
\begin{align*}
\nabla_{\boldsymbol{\lambda}}g_{i}\left(\boldsymbol{x}^{*},\boldsymbol{y}^{J*},\boldsymbol{\xi}\right) & =\partial_{\boldsymbol{\lambda}}\boldsymbol{y}^{J*}\partial_{\boldsymbol{y}}g_{i}\left(\boldsymbol{x}^{*},\boldsymbol{y}^{J*},\boldsymbol{\xi}\right),\\
\nabla_{\boldsymbol{\lambda}}h_{j}\left(\boldsymbol{y}^{J}\left(\boldsymbol{x}^{*},\boldsymbol{\lambda}^{*},\boldsymbol{\xi}\right),\boldsymbol{\xi}\right) & =\partial_{\boldsymbol{\lambda}}\boldsymbol{y}^{J*}\partial_{\boldsymbol{y}}h_{j}\left(\boldsymbol{y}^{J*},\boldsymbol{\xi}\right),
\end{align*}
where $\partial_{\boldsymbol{\lambda}}\boldsymbol{y}^{J*}=\partial_{\boldsymbol{\lambda}}\boldsymbol{y}^{J}\left(\boldsymbol{x}^{*},\boldsymbol{\lambda}^{*},\boldsymbol{\xi}\right)\in\mathbb{C}^{m\times n_{y}}$
are the derivative of the vector function $\boldsymbol{y}^{J}\left(\boldsymbol{x},\boldsymbol{\lambda},\boldsymbol{\xi}\right)$
to the vector $\boldsymbol{\lambda}$ at point $\boldsymbol{x}^{*},\boldsymbol{\lambda}^{*}$.
Therefore, after multiplying both sides of the first KKT condition
in (\ref{eq:STsconeJ}) for $\left(\boldsymbol{x},\boldsymbol{\lambda}\right)=\left(\boldsymbol{x}^{*},\boldsymbol{\lambda}^{*}\right)$
with $\partial_{\boldsymbol{\lambda}}\boldsymbol{y}^{J*}$ and taking
expectation w.r.t. $\boldsymbol{\xi}$, we have
\begin{equation}
\left\Vert \nabla_{\boldsymbol{\lambda}}f_{0}^{J}\left(\boldsymbol{x}^{*},\boldsymbol{\lambda}^{*}\right)+\sum_{i}\lambda_{i}^{*}\nabla_{\boldsymbol{\lambda}}f_{i}^{J}\left(\boldsymbol{x}^{*},\boldsymbol{\lambda}^{*}\right)\right\Vert =O\left(e\left(J\right)\right).\label{eq:PT1}
\end{equation}
In (\ref{eq:PT1}), we have used Lemma \ref{lem:keylemma1}. Combining
(\ref{eq:PT1}) and the second KKT condition in (\ref{eq:STsconeJ}),
and, we have
\begin{align}
\left\Vert \sum_{i\in\mathcal{I}_{IA}^{L}}\lambda_{i}^{*}\nabla_{\boldsymbol{\lambda}}f_{i}^{J}\left(\boldsymbol{x}^{*},\boldsymbol{\lambda}^{*}\right)-\right.\nonumber \\
\left.\sum_{i\in\mathcal{I}_{A}^{L}}\left(\widetilde{\lambda}_{i}-\lambda_{i}^{*}\right)\nabla_{\boldsymbol{\lambda}}f_{i}^{J}(\boldsymbol{x}^{*},\boldsymbol{\lambda}^{*})\right\Vert  & =O\left(e\left(J\right)\right)\label{eq:PT2}
\end{align}
where $\mathcal{I}_{A}^{L}=\left\{ i:f_{i}^{J}(\boldsymbol{x},\boldsymbol{\lambda})=0\right\} $
is the index set of the active long-term constraints, and $\mathcal{I}_{IA}^{L}=\left\{ i:f_{i}^{J}(\boldsymbol{x},\boldsymbol{\lambda})<0\right\} $
is the index set of the inactive long-term constraints. It follows
from (\ref{eq:PT2}) and the LIRC that 
\begin{align}
\lambda_{i}^{*} & =O\left(e\left(J\right)\right),\forall i\in\mathcal{I}_{IA}^{L},\nonumber \\
\widetilde{\lambda}_{i}-\lambda_{i}^{*} & =O\left(e\left(J\right)\right),\forall i\in\mathcal{I}_{A}^{L},\label{eq:PT2-2}
\end{align}
which indicates that the third KKT condition in (\ref{eq:CompS})
is satisfied up to error $O\left(e\left(J\right)\right)$. 

Similar to (\ref{eq:PT1}), it can be shown that
\begin{align}
 & \left\Vert \mathbb{E}\left[\partial_{\boldsymbol{x}}\boldsymbol{y}^{J*}\partial_{\boldsymbol{y}}g_{0}\left(\boldsymbol{x}^{*},\boldsymbol{y}^{J*},\boldsymbol{\xi}\right)\right]\right.\nonumber \\
 & +\left.\sum_{i}\lambda_{i}^{*}\mathbb{E}\left[\partial_{\boldsymbol{x}}\boldsymbol{y}^{J*}\partial_{\boldsymbol{y}}g_{i}\left(\boldsymbol{x}^{*},\boldsymbol{y}^{J*},\boldsymbol{\xi}\right)\right]\right\Vert =O\left(e\left(J\right)\right),\label{eq:PT1-1}
\end{align}
where $\partial_{\boldsymbol{x}}\boldsymbol{y}^{J*}=\partial_{\boldsymbol{x}}\boldsymbol{y}^{J}\left(\boldsymbol{x}^{*},\boldsymbol{\lambda}^{*},\boldsymbol{\xi}\right)\in\mathbb{C}^{m\times n_{y}}$
is the derivative of the vector function $\boldsymbol{y}^{J}\left(\boldsymbol{x},\boldsymbol{\lambda},\boldsymbol{\xi}\right)$
to the vector $\boldsymbol{x}$ at point $\boldsymbol{x}^{*},\boldsymbol{\lambda}^{*}$.
Applying the chain rule to the first KKT condition in (\ref{eq:LKKT}),
we have
\begin{align}
\partial_{\boldsymbol{x}}f_{0}\left(\boldsymbol{x}^{*},\Theta^{*}\right)+\sum_{i}\widetilde{\lambda}_{i}\mathbb{E}\left[\partial_{\boldsymbol{x}}\boldsymbol{y}^{J*}\partial_{\boldsymbol{y}}g_{i}\left(\boldsymbol{x}^{*},\boldsymbol{y}^{J*},\boldsymbol{\xi}\right)\right]\nonumber \\
+\sum_{i}\widetilde{\lambda}_{i}\partial_{\boldsymbol{x}}f_{i}(\boldsymbol{x}^{*},\Theta^{*})+\mathbb{E}\left[\partial_{\boldsymbol{x}}\boldsymbol{y}^{J*}\partial_{\boldsymbol{y}}g_{0}\left(\boldsymbol{x}^{*},\boldsymbol{y}^{J*},\boldsymbol{\xi}\right)\right] & =\boldsymbol{0}.\label{eq:PT3}
\end{align}
Then it follows from (\ref{eq:PT1-1}), (\ref{eq:PT2-2}) and (\ref{eq:PT3})
that
\[
\left\Vert \partial_{\boldsymbol{x}}f_{0}\left(\boldsymbol{x}^{*},\Theta^{*}\right)+\sum_{i}\lambda_{i}^{*}\partial_{\boldsymbol{x}}f_{i}(\boldsymbol{x}^{*},\Theta^{*})\right\Vert =O\left(e\left(J\right)\right),
\]
i.e., the second KKT condition in (\ref{eq:CompS}) is satisfied up
to error $O\left(e\left(J\right)\right)$. 

Finally, it follows directly from (\ref{eq:STsconeJ}) that the first
KKT condition in (\ref{eq:CompS}) is satisfied up to error $O\left(e\left(J\right)\right)$.
This completes the proof.

\subsubsection{Proof of Lemma \ref{lem:Property-surrogate} \label{subsec:Proof-of-Lemma-surr}}

From Assumptions \ref{asm:convP} and \ref{asm:shortAlg}, $\partial_{\boldsymbol{x}}^{T}g_{i}^{\bar{c}}\left(\boldsymbol{x}^{t},\boldsymbol{y}_{j}^{t},\boldsymbol{\xi}_{j}^{t}\right)$,
$\mathbf{f}_{x,i}^{t-1}$, $\mathbf{f}_{y,i}^{t}$ and $\mathbf{f}_{\lambda,i}^{t}$
are bounded. Therefore, Results 1) - 3) in Lemma \ref{lem:Property-surrogate}
follow directly from the expression of the structured surrogate function
in (\ref{eq:SSF}). Result 4) is ture because $\lim_{t\rightarrow\infty}f_{i}^{t}-f_{i}^{J}(\boldsymbol{x}^{t},\boldsymbol{\lambda}^{t})=0$,
$\lim_{t\rightarrow\infty}\left\Vert \mathbf{f}_{x,i}^{t}+\mathbf{f}_{y,i}^{t}-\nabla_{\boldsymbol{x}}f_{i}^{J}(\boldsymbol{x}^{t},\boldsymbol{\lambda}^{t})\right\Vert =0$
and $\lim_{t\rightarrow\infty}\left\Vert \mathbf{f}_{\lambda,i}^{t}-\nabla_{\boldsymbol{\lambda}}f_{i}^{J}(\boldsymbol{x}^{t},\boldsymbol{\lambda}^{t})\right\Vert =0$,
which is a consequence of the chain rule and (\cite{Ruszczyski_MP80_SPthem},
Lemma 1). The proof is similar to that of (\cite{Yang_TSP2016_SSCA},
Lemma 1) and the details are omitted for conciseness. Finally, Result
5) is a consequence of Result 4).


\end{document}